\documentclass[a4paper,12pt]{article}
\usepackage{ulem}
\usepackage{dmvn}
\usepackage{graphicx}
\usepackage[matrix,arrow,curve]{xy}
\usepackage{lscape}

\newtheorem{opr}{Определение}
\newtheorem{zamech}{Замечание}
\newtheorem{utv}{Утверждение}
\newtheorem{prim}{Пример}
\newtheorem{sled}{Следствие}
\newtheorem{predl}{Предложение}

\begin{document}
\author{Andrey Mikhovich\footnote{Corresponding author.
E-mail address: mikhandr@mail.ru}}
\title{Homotopy of profinite groups.}
\date{28.04.12}
\maketitle
\begin{abstract}
We study simplicial profinite groups with a view towards applications in profinite combinatorial group
theory. This approach provides a natural framework to the concept of pro-$\mathfrak{C}$-presentation of a pro-$\mathfrak{C}$-group $G$
as a 1-truncation of its free simplicial pro-$\mathfrak{C}$-resolution. The category of simplicial pro-$\mathfrak{C}$-groups has closed
simplicial model category structure. This yields a possibility to define some old and new derived functors as left Quillen derived
functors from this simplicial model category. When $\mathfrak{C}$ is a L-groups we construct free simplical pro-L-resolution
functorially. We introduce settings of $\Delta$-adic and Zassenhaus filtrations for free simplical pro-p-resolutions and derive some
calculations for pro-p-groups. The usage of pro-p-Curtis-Rector spectral sequences sheds homotopical light on Golod-Shafarevitch type
results.
\end{abstract}

\section{Введение.}
Работа посвящена методам симплициальной теории про-$\mathcal{C}$-групп с акцентом на возможные приложения в комбинаторной теории.
Отправной точкой стала статья   \cite{Sta} , в которой построена (довольно запутанно) некоторая спектральная последовательность,
связанная с копредставлением дискретной группы, а также исследована задача существования эпиморфизма некоторых дискретных групп на
свободные группы конечного ранга.

Систематический  анализ вопроса в про-$\mathcal{C}$-случае (мы сразу ограничились многообразием $\mathcal{C}$ конечных групп «полным» в
смысле \cite{AM}, так как это впоследствии позволяет строить свободные симплициальные резольвенты функториально) привёл к выводу, что для
построения подобного типа спектральных последовательностей необходимо отталкиваться от про-$\mathcal{C}$-симплициальной резольвенты
про-$\mathcal{C}$-группы, построенной по заданному про-$\mathcal{C}$-копредставлению.
Задача построения реализована в работе с помощью модификации метода, введенного Андре в \cite{And}. Специфика про-$\mathcal{C}$-ситуации
проявляется в удобстве использования пунктированных проконечных пространств, заменяющих образующие в свободных дискретных группах.  В
определении заложено существовании $CW$-базиса, т.к. в случае про-$\mathcal{C}$-групп, в отличие от дискретной ситуации (есть результат
Кана), его наличие, вообще говоря,  не гарантировано.

Имея на руках симплициальные про-$\mathcal{C}$-резольвенты, построенные методом Андре, интересно выяснить их общее место в категории
симплициальных про-$\mathcal{C}$-групп. Этот вопрос привел к рассмотрению данной категории с позиций гомотопической алгебры. Доказано,
что категория симплициальных про-$\mathcal{C}$-групп является симплициальной модельной категорией со «стандартным», а не как например у
Мореля- гомологическая локализация, набором морфизмов. Структура симплициальной модельной категории обеспечивает нас  критерием
существования производных функторов по Квиллену, что создает базис для построения как абелевой, так и неабелевой гомологической алгебры.
А свободные симплициальные резольвенты представляют собой кофибрантные замены, на которых и вычисляются производные функторы. Рассмотрены
примеры некоторых классических и новых, в том числе, неабелевых производных функторов на  категории симплициальных
про-$\mathcal{C}$-групп. С помощью идеи, содержащейся в работе M.Artin, B.Mazur, On the Van Kampen theorem , Topology 5(1966) ,
используя спектральную последовательность Квиллена,   доказано например, что обычные гомологии про-$p$-группы с коэффициентами  в
$\mathbb{F}_p$- суть производные функтора факторизации по подгруппе Фраттини.

   В третьей  части работы рассмотрены  фильтрации свободных симплициальных резольвент про-$p$-групп и соответствующих групповых колец.
   Оказывается, вопрос о наличии сюръективного гомоморфизма в свободную про-р-группу заданного ранга имеет гомотопическую природу и
   состоит в нахождении ретракции из свободной симплициальной резольвенты, построенной по заданному копредставлению, в «постоянную»
   свободную про-$p$-группу заданного ранга. При этом использование информации, возникающей из фильтрации, приводит  к конкретным
   примерам и результатам. Задачи такого типа имеют тесную связь с гипотезой Гринберга, например  K.Wingberg, Free pro-p-extensions of
   number fields, 2005  (и цитируемая в этой работе литература).

Четвертая  часть работы имеет феноменологический характер и объясняет связь между  свободными симплициальными резольвентами,
резольвентами скрещенных модулей и, наконец, свободными резольвентами в категории компактных $G$– модулей.Самостоятельный интерес
представляет чисто симплициальное доказательство того, что $\pi_2 (P)$ – вторая гомотопическая группа скрещенного модуля копредставления
некоторой про-$p$-группы $G$ (в дискретной ситуации это $\pi_(\mathfrak{С}_P)$ – вторая гомотопическая группа стандартного
$CW$–комплекса, полученная из копредставления $P$ группы $G$ ) изоморфна первой гомотопической группе комплекса Мура второго шага
построения резольвенты методом «pas a pas». Представленные в литературе доказательства  используют так называемые «identity  sequences»
Pride S.J. Identities  among  relations  of  group  presentations 1991, которые не могут быть использованы в про-$\mathcal{C}$-случае.

  Про-$p$-группы естественно возникают как про-$p$-пополнения дискретных групп и дают возможность взглянуть на естественные структуры
  по-новому, словами Дж.Адамса, сквозь различные «$p$-адические очки». Однако такая  идея натыкается на сложные нерешенные проблемы в
  самой теории про-$p$-групп. В частности, мы обсуждаем старую открытую проблему о про-$p$-группах с одним соотношением, поставленную
  Серром на семинаре Бурбаки. Многократно подчеркивается значение и место модулей соотношений, как ключевого и наиболее важного с
  практической точки зрения объекта. Знание модуля соотношений дает теоретическую возможность с помощью сдвига размерности вычислять все
  (ко)гомологии через первые и вторые. Зачастую, как в случае групп с одним соотношением, модуль соотношений определяет и весь
  гомотопический тип.

В последней части работы исследуются спектральные последовательности, возникающие из фильтраций на свободной симплициальной резольвенте
про-$p$-группы $G$ . Оказывается, что это и есть системный  взгляд на работу Столлингса, а по сути исследование нестабильных спектральных
последовательностей Адамса, которые ранее встречались в работах Кертиса, Ректора, Боусфилда, позже Грюненфельдера. Выведены
про-$p$-варианты основных теорем, исследована сходимость. Замечено, что вопрос Тима Портера о природе неравенства
Голода-Шаферевича-Винберга-Коха тоже имеет гомотопическую природу и может быть выражен как характеристика скорости сходимости полученной
спектральной последовательности.

Автор выражает благодарность Владимиру Барановскому за совместные беседы, которые не только давали пищу для дальнейших размышлений, но и
оказывали влияние за рамками математики. Я признателен Тиму Портеру и Роману Михайлову за ранние версии \cite{P1} и \cite{Mikh}, которые
оказались полезными при написании этой работы. Своему интересу к про-р-группам автор обязан Олегу Владимировичу Мельникову.

Спасибо моей жене Елене за моральную поддержку, уют и помощь в работе с LaTeX.

\section{Симплициальные структуры.}
Прежде всего изложим хорошо известные построения, чтобы установить некоторые обозначения и терминологию, которые нам понадобятся в
дальнейшем.

Категория $\Delta:  Ob(\Delta) = \{[n] = \{0,1,...n\}\}, Mor(\Delta) = \{\mu: [n]\to [m] | \, \mu i \le \mu j, i \le j \}$

Пусть $\mathcal{C}$ - произвольная категория. Контравариантный функтор $S: \Delta \to \mathcal{C}$ будет называться симплициальным
объектом в $\mathcal{C}$. В нашем случае это значит, что $S$ ставит в соответствие каждому неотрицательному целому числу $n$ (каждому
объекту из $\Delta$) объект $S_n$ из $\mathcal{C}$ и каждому монотонному отображению $\mu: [n] \to [m]$ морфизм $\mu^*=S(\mu): S_m \to
S_n$ из $\mathcal{C}$, причем $S(1)=1$ и $S(\mu\nu)=S(\nu) S(\mu).$ Под симплициальным множеством будет пониматься симплициальный объект
в категории множеств SS; под симплициальной про-$\mathfrak{C}$-группой будет пониматься симплициальный объект в категории
про-$\mathfrak{C}$-групп $sGr^{pro-\mathfrak{C}}$, под редуцированным симплициальным множеством будем понимать симплициальное множество,
у которого $SS_0$ состоит из одного элемента $\{*\}$ и т.д.

\textbf{Класс $\mathfrak{C}$ будем предполагать "полным" в смысле \cite{AM} или, что почти тоже самое, что $\mathfrak{C}$ является
многообразием конечных групп, замкнутым относительно расширений \cite{ZR}. В частности, аналогично \cite{AM}, нас будут в основном
интересовать случаи, когда $\mathfrak{C}$ либо класс всех конечных групп, либо класс всех конечных $L$-групп, т.е. конечных групп,
порядки которых есть произведения степеней из набора простых чисел $L$.}

Пусть $[n]$ - это упорядоченное множество, $[n]={0<1<...<n}$. Определим следующие отображения. Во-первых, инъективное монотонное
отображение $\delta_i^n: [n-1]\to [n]$, заданное

$$
\delta_i^n(x)=\begin{cases}
x,&\text{если $x<i$;}\\
x+1,&\text{если $x \ge i$,}
\end{cases}
$$
для $0 \le i \le n \ne 0.$ Для удобства будем изображать и писать эти отображения, опуская верхние индексы:

$$
\xymatrix{
[0]\ar@<1ex>[r]^{\delta_0}\ar@<-1ex>[r]^{\delta_1}&[1]\ar@<2ex>[r]^{\delta_0}\ar@<0ex>[r]^{\delta_1}\ar@<-2ex>[r]^{\delta_2}&[2]\ar@<3ex>[r]^{\delta_0}\ar@<1ex>[r]^{\delta_1}\ar@<-1ex>[r]^{\delta_2}\ar@<-3ex>[r]^{\delta_3}&[3]
...}$$
Возрастающее сюръективное монотонное отображение $\al_i^n: [n+1]\to [n]$, заданное
$$
\al_i^n(x)=\begin{cases}
x,&\text{если $x \le i$;}\\
x-1,&\text{если $x > i$,}
\end{cases}
$$
для $0 \le i \le n.$
$$
\xymatrix{
[0]&[1]\ar[l]_{\al_0}&[2]\ar@<1ex>[l]_{\al_0}\ar@<-1ex>[l]_{\al_1}&[3]\ar@<-2ex>[l]_{\al_0}\ar@<0ex>[l]_{\al_1}\ar@<2ex>[l]_{\al_2}...}
$$

Немедленно из определений вытекают следующие равенства
\begin{eqnarray}
\delta_j \delta_i= \delta_i \delta_{j-1}, i < j,\\
\al_j \al_i = \al_i \al_{j+1}, i \le j,\\
\al_j \delta_i = \delta_i \al_{j-1}, i < j,\\
\al_j \delta_j = \al_j \delta_{j+1}=id,\\
\al_j \delta_i = \delta_{i-1} \al_j, i > j+1.
\label{1}
\end{eqnarray}

\begin{lemma} ~\cite[гл.~2.1]{GZ}
Каждая монотонная функция $\mu: [n] \to [m]$ представима единственным образом в виде произведения
\begin{eqnarray}
\mu = \delta^{i_1} ... \delta^{i_s} \al^{j_1} ...\al^{j_t},
\mbox{ где } m \ge i_1 >...>i_s \ge 0, \, 0 \le j_1 < ... < j_t < n, \, n-t+s = m.\label{45}
\end{eqnarray}
\end{lemma}

Из этой леммы немедленно вытекают два следствия - см.~\cite[гл.~2.1]{GZ}.

\begin{sled}
Категорию $\Delta$ можно отождествить с категорией $\Delta'$, порожденной объектами $[n]$ и морфизмами $\delta, \al$, подчиненными
соотношениям (1)-(5).
\end{sled}

\begin{sled}
(Альтернативное определение симплициального объекта)
Симплициальный объект F в категории $\mathcal{C}$ - это семейство $\{S_n\}$ объектов из $\mathcal{C}$ вместе с двумя семействами
морфизмов из $\mathcal{C}$
$$d_i: S_n \to S_{n-1}, s_i: S_n \to S_{n+1}, \, i=0,...,n$$
(и n>0 в случае $d_i$), которые удовлетворяют равенствам

\begin{eqnarray}
d_i d_j = d_{j-1} d_i, \,i < j,\\
s_i s_j = s_{j+1} s_i, \,i \le j,\\
d_i s_j = s_{j-1} d_i, \,i < j,\\
d_j s_j = d_{j+1} s_j = id\\
d_i s_j = s_j d_{i-1}, \, i>j+1.\label{2}
\end{eqnarray}
\end{sled}
Мы назовем $d_i$ и $s_j$ соответственно $i$-м граничным оператором и $j$-м оператором вырождения симплициального объекта $S$. Из (7)-(8)
вытекают равенства
\begin{eqnarray}
d_i d_j = d_j d_{i+1}, \,i \ge j,\\
s_i s_j = s_j s_{i-1}, \,i > j.\label{3}
\end{eqnarray}

\begin{prim}
Симплициальное множество $\Delta[n]$.

$\Delta[n]_m = Mor_{\Delta}([m],[n])$

$d_i^m: \Delta[n]_m \to \Delta[n]_{m-1} $ определяется, как композиция:
$[m-1] \xrightarrow{\al_i} [m] \xrightarrow{\gamma} [n]$

$s_i^m: \Delta[n]_m \to \Delta[n]_{m+1} $ определяется, как композиция:
$[m+1] \xrightarrow{\delta_i} [m] \xrightarrow{\gamma} [n]$
Ясно, что у симплициального множества $\Delta[n]$ есть только один невырожденный симплекс в размерности $n$, все остальные элементы
$\Delta[n]$ получаются либо как грани, либо как вырождения этого $n$-мерного симплекса.  $\Lambda_n^k$ - симплициальное подмножество
$\Delta[n]$, порожденное всеми $(n-1)$-мерными гранями $\Delta[n]$, за исключением $k$-ой.
\end{prim}

\begin{opr}
Cимплициальное множество $G$ обладает свойством Кана и называется "фибрантным", если любой морфизм $\varphi: \Lambda_n^k \to G$
продолжается до морфизма $\psi: \Delta [n] \to G$.
\end{opr}
Поскольку морфизм $\varphi$ полностью определяется элементами $\varphi(d^n_i) \in G_{n-1}, \, i \ne k$, а морфизм $\psi$-элементом
$\psi(Id_{[n]}) \in G_n$, свойство Кана равносильно тому, что если заданы элементы $x_0,...,x_{k-1},x_{k+1},...,x_n \in G_{n-1}, \,
d^{n-1}_i x_j = d^{n-1}_{j-1} x_i, \, 0 \le i \le j \le n, \, ci,j \ne k$, то существует $x \in G_n: d^n_i (x)= x_i$.

\begin{predl}
Cимплициальная группа является множеством Кана.
\end{predl}

\textbf{Доказательство.} ~\cite[гл.~4]{May}.

\begin{opr}
Комплексом Мура $(NG, \partial)$ симплициальной группы $G$ назовем цепной комплекс:
$$
NG_n = \bigcap_{i=0}^{n-1} \Ker d_i^n, \, \partial_n: NG_n \to NG_{n-1},$$ $ \partial_n$ - это ограничение $d_n^n$ на $NG_n$, $\pi_n(G)$
- $n$-ая гомотопическая группа определяется как n-ые гомологии комплекса Мура
$
\pi_n(G) := H_n (NG, \partial) =\frac{\bigcap\limits_{i=0}^n \Ker d_i^n} {d_{n+1}^{n+1} (\bigcap\limits_{i=0}^n \Ker d_i^{n+1})}
$
\end{opr}

В работе \cite{Qui4} Квиллен свободно использует технологию присоединения новых элементов к симплициальной проконечной группе ("cell
attachment procedure"), никак не поясняя её природу. Дадим этому процессу определённость и покажем, что такая процедура может быть
использована для построения свободной про-$\mathfrak{C}$-резольвенты про-$\mathfrak{C}$-группы $G$, как неотъемлемая составляющая
неабелевого про-$\mathfrak{C}$-аналога метода "pas-a-pas" \cite{And}.
Пусть $G$ - симплициальная про-$\mathfrak{C}$-группа и $\Omega$ - замкнутое пунктированное (отмеченная точка - это единица группы)
подпространство в $N G_{k-1}$.
Мы хотим построить симплициальную про-$\mathfrak{C}$-группу $F = G(\Omega)$ с мономорфизмом
$\beta : G \to F$
таким, что гомоморфизм
$\pi_{k-1}(\beta): \pi_{k-1}(G) \to \pi_{k-1} (F)$ равен нулю на элементах из $ \Omega.$ Это будет сделано в теореме 1.2 с помощью
добавления новых элементов $x_{\la}$ к образу $\partial_k (NG_k)$ таким образом, чтобы $\beta (x_{\la}) \in \partial NF_k.$
 В последующем изложении $\Omega$ будет возникать из замкнутой пунктированной системы образующих $\Omega'$ про-$\mathfrak{C}$-группы
 $\pi_{k-1}(G)$. Пусть $p: N G_{k-1} \twoheadrightarrow \pi_{k-1}(G)$ - каноническая сюръекция. По ~\cite[теор.~1]{Se1} у $p$ существует
 непрерывное сечение $\tau : \pi_{k-1}(G) \to N G_{k-1}$. Т.к. $\Omega'$ замкнуто в $\pi_{k-1}(G)$, то оно бикомпактно \cite[\S
 13,(B)]{Pon}, тогда $\tau (\Omega')$ бикомпактно \cite[\S 13,(D)]{Pon}, а следовательно \cite[\S 13,(C)]{Pon} $\tau (\Omega')$-
 замкнутое подпространство в $ N G_{k-1}.$
 В про-р-случае выбор такой системы образующих может быть сделан ещё более каноническим способом, а именно: выберем замкнутую
 пунктированную систему $\Omega''$ образующих топологического векторного пространства $\pi_{k-1}(G)/\Phi(\pi_{k-1}(G))$ над полем
 $\mathbb{F}_p$ (k > 1, $\Phi$ - подгруппа Фраттини). Т.к. у канонической сюръекции $\pi_{k-1}(G) \twoheadrightarrow
 \pi_{k-1}(G)/\Phi(\pi_{k-1}(G))$ существует непрерывное сечение $\sigma$, то $\Omega' = \sigma (\Omega'')$ - замкнутая система
 образующих $\pi_{k-1}(G)$. Примем $\Omega = \tau \sigma(\Omega'')$.  Подробности относительно свободных про-$\mathfrak{C}$-групп можно
 найти в \cite{ZR} и \cite{Gl1}.

  Ключевым результатом этой части является следующая теорема, которая является аналогом результата Андрэ \cite[Предл. 6.1]{And}

\begin{theorem}
Пусть $G$ - симплициальная про-$\mathfrak{C}$-группа, задано $k>1$ и замкнутое пунктированное подпространство $\Omega\subseteq
\pi_{k-1}(G) $,  $\Omega = \{x_{\la} \in \pi_{k-1}(G) | \, \la \in \Lambda\}$.
Тогда правила 1)-3) корректно определяют структуру симплициальной про-$\mathfrak{C}$-группы $(F_n, s_i^n, d_i^n)$:

1) $F_k = G_k * F_{id[k]}(\Omega),$
$ F_n=G_n *F(\cup_{t:[n]\twoheadrightarrow [k]} \Omega_t), \quad \Omega_t \cong \Omega$

$F_n=\begin{cases}
G_n,&n < k,\\
G_k*F(\Omega),& n=k,\\
G_n*F(\cup_{t:[n]\twoheadrightarrow [k]}\Omega_t),\quad \Omega_t \cong \Omega
\end{cases}
$

Если $Y$- пунктированное проконечное пространство, то мы обозначаем $F(Y)$ свободную про-$\mathfrak{C}$-группу, порожденную $Y$.

2) Для $0 \le i \le n$, гомоморфизм про-$\mathfrak{C}$-групп $s_i^n: F_n \to F_{n+1}$ получен из гомоморфизма $s_i^n: G_n \to G_{n+1}$ на
новых образующих по правилу
$$s_i^n (y_{\la,t})=y_{\la,u} \text{, где } u=t\al_i^n, \, t:[n] \to [k], \, y_{\la,t}\in \Omega_t $$

3) Для $0 \le i \le n \ne 0,$ гомоморфизм про-$\mathfrak{C}$-групп $d_i^n: F_n \to F_{n-1}$ получен из $d_i^n: G_n \to G_{n-1}$ с помощью
задания на новых образующих по правилу
$$
d_i^n(y_{\la,t})=\begin{cases}
y_{\la,u},&\text{если отображение $u=t\delta_i^n$ является сюръективным,}\\
\tilde t'(\vartheta_{\la}),&\text{если $u=\delta_k^k t'$,}\\
1,&\text{если $u=\delta_j^k t'$ с $j \ne k,$}
\end{cases}
$$
$t': [n-1] \to [k-1]$ морфизм в $\Delta$, которому соответствует единственный гомоморфизм про-$\mathfrak{C}$-групп $\tilde t': G_{k-1}
\to G_{n-1}.$

\end{theorem}

\begin{opr}
Аугментированная симплициальная про-$\mathfrak{C}$-группа $(F,\varepsilon,G)$ называется асферичной, если $\pi_n(F)=0$ при $n>0$ и
$\pi_0(F)=G$
\end{opr}

\begin{opr}
Свободной симплициальной резольвентой про-$\mathfrak{C}$-группы $G$ называется асферичная про-$\mathfrak{C}$-группа $X$ с $\pi_0(X)=G,$ у
которой имеется $CW$-базис \cite{Kan2}, т.е.

(a) $X_n$ свободно порождается проконечным пространством $\mathcal{F}_n$

(b) $\mathcal{F}_n$ устойчивы относительно операторов вырождения, т.е. если $\sigma \in \mathcal{F}_n,$ то $s_i(\sigma) \in
\mathcal{F}_{n+1}, \, 0 \le i \le n, \, \mathcal{F}_{n+1} = \cup_{i=1}^n s_i(\mathcal{F}_n) \uplus Y_{n+1},$ где $Y_{n+1}$ - замкнутое
подпространство в $X_{n+1}.$
\end{opr}

\begin{theorem}

У каждой про-$\mathfrak{C}$-группы $G$ есть свободная симплициальная  про-$\mathfrak{C}$-резольвента.
\end{theorem}
\textbf{Доказательство.}

Пусть $G$ - группа. Опишем нулевой шаг конструкции. Он состоит из выбора свободной про-$\mathfrak{C}$-группы $F$ и сюръекции $g: F \to
G$, что дает изоморфизм $F/\Ker g := G$, как групп. Образуем постоянную симплициальную про-$\mathfrak{C}$-группу $F^{(0)}$, для которой
для каждой степени $n$ выполнено $F_n = F$ и $d_i^n = id = s_j^n$ для всех $i,j$, а $\pi_0(F^{(0)}) = G.$ Выберем замкнутое
подпространство $\Omega^{(0)}$ нормальных образующих для $N=\Ker(F \xrightarrow{g}G)$ и рассмотрим симплициальную
про-$\mathfrak{C}$-группу, в которой $F_1^{(1)}: = F(X_0 \cup \Omega^{(0)})$ и $F_n^{(1)}$ для $n > 1$ будет свободной
про-$\mathfrak{C}$-группой на пространстве вырожденных образующих, описанных выше. Эту симплициальную про-$\mathfrak{C}$-группу обозначим
$F^{(1)}$ и будем называть 1-скелетом симплициальной резольвенты про-$\mathfrak{C}$-группы $G$.
Последующие шаги зависят от выбора $\Omega^0,\Omega^1,\Omega^2,...,\Omega^k,...$ Пусть $F^{(k)}$ - симплициальная
про-$\mathfrak{C}$-группа, полученная после $k$-ого шага, $k$-ого скелета резольвенты. Пространство $\Omega^k$ образовано с помощью
элементов, соответствующих образующим абелевой про-$\mathfrak{C}$-группы $\pi_k(F^{(k)})$.

Имеем включения симплициальных про-$\mathfrak{C}$-групп:
$F^{(0)}\subseteq F^{(1)}\subseteq...\subseteq F^{(k-1)}\subseteq F^{(k)}\subseteq...,$
и при переходе к индуктивному пределу, мы получаем ацикличную свободную симплициальную про-$\mathfrak{C}$-группу $F$, $\pi_0(F) = G$.
$\boxtimes$

Чтобы сделать предыдущие построения еще более наглядными, построим симплициальную резольвенту вплоть до размерности два. Пусть $X_0$ -
некоторое пунктированное пространство образующих $G$. Зафиксируем каноническую сюръекцию $\varepsilon : F(X_0)\twoheadrightarrow G$.
1-скелет свободной симплициальной резольвенты про-$\mathfrak{C}$-группы $G$ строится путем добавления новых образующих $Y_1$, находящихся
в 1-1 соответствии с нормальными образующими $\Ker  \varepsilon$ (как топологического нормального делителя), тогда $F_1^1 = F(X_0 \cup
Y_1) $, и получаем:

$$
\xymatrix{{F(X_0 \cup Y_1)} \ar@<2ex>[rr]^{d_1^1} \ar@<0ex>[rr]^{d_0^1} && {F(X_0)} \ar@<2ex>[ll]_{s_0^0} \ar[rr]^{\varepsilon} && G
},$$

где $F(X_0) \xrightarrow{\varepsilon} G$ - гомоморфизм копредставления, а $s_0,d_0^1,d_1^1$ заданы формулами:
$$d_1^1(y_i) = b_i \in \Ker \varepsilon, \, d_0^1 (y_i) = 1,\, s_0(x_0) = x_0, \, x_0 \in X_0$$
Заметим, что $\langle X_0 | \, d_1^1 (Y_1)\rangle$ - копредставление $G$ в обычном смысле. Теперь 1-скелет $F^{(1)}$ выглядит так:

$\xymatrix{{F^{(1)}: ... F(X_0 \cup s_0(Y_1) \cup s_1 (Y_1))} \ar@<3ex>[rr]^(0.6){d_0^2,d_1^2,d_2^2} \ar@<2ex>[rr] \ar@<1ex>[rr] &&
{F(X_0 \cup Y_1)} \ar@<3ex>[ll]^(0.4){s_1^1,s_0^1} \ar@<2ex>[ll] \ar@<3ex>[rr]^{d_1^1,d_0^1} \ar@<1ex>[rr] && {F(X_0)
\ar@<1ex>[ll]^{s_0^0} \ar[r]}& G
}$
\begin{zamech}
Метод "pas a pas" может быть использован для того, чтобы построить для любой симплициальной про-$\mathfrak{C}$-группы $G$ свободную
симплициальную про-$\mathfrak{C}$-группу $F$ и эпиморфизм $F\twoheadrightarrow G,$ который индуцирует изоморфизмы гомотопических групп
\cite{MP}.
\end{zamech}

\begin{opr} \cite[гл.~6]{May}
Пусть $f$ и $g$ - два гомоморфизма симплициальных про-$\mathfrak{C}$-групп $G \to H$. Тогда они называются гомотопными и обозначаются
$f\simeq g$, если существует набор гомоморфизмов $h_i: G_q \to H_{q+1}, 0 \le i \le q,$ которые удовлетворяют следующим тождествам:
\begin{eqnarray}
d_0 h_0 = f,\\
d_{q+1} h_q = g,\\
d_i h_j = h_{j-1} d_i, i < j,\\
d_{j+1} h_j = d_{j+1} h_{j+1},\\
d_i h_j = h_j d_{i-1}, i > j+1,\\
s_i h_j = h_{j+1} s_i, i \le j.\\
s_i h_j = h_j s_{i-1}, i >j.\label{4}
\end{eqnarray}
\end{opr}

\begin{opr}Пополненной симплициальной про-$\mathfrak{C}$-группой называется симплициальная про-$\mathfrak{C}$-группа $(G_*,\varepsilon,
G)$ и $\varepsilon: G_0 \to G: \varepsilon d_0^1 = \varepsilon d_1^1$. Она называется стягиваемой слева, если существуют $h_n: G_n \to
G_{n+1}, n \ge 0$ и $h: G \to G_0:$
$\begin{cases}
\varepsilon h = 1,\\
d_0^{n+1} h_n = 1,\\
d_1^1 h_0 = h \varepsilon,\\
d_i^{n+1} h_n = h_{n-1} d_{i-1}^n, i,n \ge 1
\end{cases}$

и называется стягиваемой справа, если существуют непрерывные гомоморфизмы $h_n$ и $h$: $\begin{cases}
\varepsilon h = 1,\\
d_{n+1}^{n+1} h_n = 1, n \ge 0\\
d_0^1 h_0 = h \varepsilon,\\
d_i^{n+1} h_n = h_{n-1} d_{i}^n, n \ge 1, i \le n
\end{cases}$
\end{opr}

Из определений немедленно следует:

\begin{utv}
Стягиваемая справа (слева) симплициальная про-$\mathfrak{C}$-группа является асферичной.
\end{utv}

\begin{opr}
Отображение $f: G \twoheadrightarrow H$  симплициальных про-$\mathfrak{C}$-групп называется расслоением, если

1) $f$ - симплициальный гомоморфизм (т.е. перестановочен с операторами граней и вырождения)

2) $\forall n \ge 0$ гомоморфизмы $f_n: G_n \to H_n$ являются сюръективными.
\end{opr}

\begin{theorem} \cite{Ina}
Если $f: G \twoheadrightarrow H$ - расслоение, то имеет место точная последовательность гомотопических групп
$... \to \pi_{n+1}(H) \to \pi_n(\Ker f) \to \pi_n(G) \to \pi_n(H) \to ...$
\end{theorem}

В доказательстве гомотопической эквивалентности свободных симплициальных резольвент Кейне \cite{Keu} использует несколько иное
определение асферичности симплициального объекта. Напомним его построение.

Пусть $A \in SGr$, тогда объекты $Z_nA \in A, n \ge 1$ определяются так:
$$Z_nA = \{(a_0, ..., a_{n+1}) \in A_n^{n+2} | d_i a_j = d_{j-1} a_i, i < j\},$$
как подобъект в $A_n^{n+2}$. Для каждого $n$ существует гомоморфизм $$d = d^{(n)} : A_{n+1} \to Z_n A, da = (d_0 a, ...,d_{n+1} a), a \in
A_{n+1}.$$
Тогда аугментированный симплициальный объект $A \in SGr$ будет называться асферичным, если для каждого $n \ge -1$ гомоморфизм $d^{(n)}$
является сюръекцией. Оказывается, что определение асферичной симплициальной про-$\mathfrak{C}$-группы, данное нами выше, как
симплициальной про-$\mathfrak{C}$-группы с тривиальными гомологиями комплекса Мура, эквивалентно определению, которое использует Кейне.
Действительно, пусть $a_0,...,a_n,\widehat{a_{n+1}} \in A_n^{(n+2)}: d_i^n a_j = d_{j-1}^n a_i$ для $i < j$
Тогда из фибрантности симплициальной группы следует, что существует $y: d_i^{n+1}(y)=a, 0 \le i \le n$. Рассмотрим элемент
$d_{n+1}^{n+1}(y)$, он не обязательно равен $a_{n+1}$. Рассмотрим элемент $d_{n+1}(y)^{-1} a_{n+1}$. Оказывается, что $d_i(d_{n+1}
(y)^{-1} a_{n+1})=1$ действительно, $d_i(d_{n+1}(y)^{-1} a_{n+1}) = d_i d_{n+1}(y)^{-1} \cdot d_i(a_{n+1}) = d_n d_i(y)^{-1} \cdot d_n
a_i = 1.$
Тогда из асферичности комплекса Мура существует элемент $z \in A_{n+1}: d_i^{n+1}(z)=1, 0 \le i \le n+1$ и
$d_{n+1}^{n+1}(z)=d_{n+1}^{n+1}(y)^{-1} a_{n+1}$
Положим $x=y \cdot z,$ тогда $d_i^{n+1}(y \cdot z) = a_i$ и $d_{n+1}^{n+1}(y \cdot z) = d_{n+1}^{n+1}(y) d_{n+1}^{n+1}(y)^{-1} \cdot
a_{n+1} = a_{n+1}.$
Утверждение о том, что из определения асферичности, которое использует Кейне, следует наше очевидно, т.к. можно применить свойство
сюръективности $d^{(n)}$ к элементу $(1,...,1,g),$ где $g \in NG_n$. Теперь можно обратиться к аналогу результата Кейне о гомотопической
эквивалентности симплициальных резольвент заданной группы, доказательства которых получаются компиляцией в про-$\mathfrak{C}$-случай:

\begin{theorem}~\cite{Keu}
Пусть $A \in asA$ - свободная симплициальная про-$\mathfrak{C}$-группа с $CW$-базисом, $B$ - асферическая и пусть $\al: A_{-1} \to
B_{-1}$ - гомоморфизм.

1) Тогда существует $\gamma : A \to B :\gamma_{-1}=\al$

2) Если $\gamma'$ - тоже расширяет $\al$, тогда $\gamma$ и $\gamma'$ гомотопны. $\boxtimes$\label{22}
\end{theorem}
В случае конечнопорождённой про-р-группы можно утверждать несколько большее. Дело в том, что элементы про-р-группы порождают
(топологически) группу $G$ тогда и только тогда, когда они дают базис векторного $\mathbb{F}_p$-пространства $X$, $G/Ф(G) \cong
G/\overline{G^p}[G,G] \cong \langle X\rangle_{\mathbb{F}_p}$ (см. ~\cite{ZR}). Тогда мы можем отождествить две системы образующих при
помощи линейного преобразования $\al \in G L_{rank(X)}(\mathbb{F}_p)$, которое можно продолжить до изоморфизма $\tilde \al: F(X_1) \cong
F(X_2)$, т.к. ввиду вышесказанного $|X_1|=|X_2|$. Похожая ситуация имеет место для локальных колец в последней главе~\cite{T}.

Пусть $G$ - симплициальная про-$\mathfrak{C}$-группа. Определяем симплициальное проконечное пространство $WG$
~\cite[опр.~6.14]{Cur71},\cite{GJ}, полагая:
$$W_n G = G_n \times G_{n-1} \times ... \times G_0, n \ge 0$$
и операторами граней и вырождения
$$d_i (g_n,...g_0)=(d_i g_n, d_{i-1} g_{n-1},...,d_0 g_{n-i} \cdot g_{n-i-1}, g_{n-i-2},...,g_0), i < n$$
$$d_n (g_n,...,g_0) = (d_n g_n, ...,d_1 g_1)$$
$$s_i (g_n,...g_0)=(s_i g_n, s_{i-1} g_{n-1},..., s_0 g_{n-i}, 1, g_{n-i-1}, g_{n-i-2},...,g_0).$$
Симплициальное пространство $WG$ имеет естественную структуру $G$-множества, именно, $G$ действует свободно слева по правилу:
$$g(g_n,...,g_0) = (g g_n,g_{n-1},...,g_0), g \in G_n .$$
Тогда $\overline{W}G$ - это факторпространство $WG$ по левому $G$-действию; обозначим отображение факторизации (классифицирующее
$G$-расслоение):
$$q = q_G: WG \to \overline{W}G$$
Конструкция $\overline{W}G$ определяет функтор:
$\overline{W}: SGr^{pro-\mathfrak{C}} \to SS_0^{profinite}$ (категория $SS_0^{profinite}$ - это категория симплициальных проконечных
пространств, у которых нулевая компонента состоит из одного элемента),
который будем называть функтором классифицирующего пространства категории про-$\mathfrak{C}$-групп. Ясно \cite[IV, \S 21]{May}, что
$\overline{W}G$ может быть описано, как $\overline{W}_0 G = [\quad], \overline{W}_n G=(e) \times G_{n-1} \times ... \times G_0, n > 0.$
Операторы граней и вырождения определяются, как $$d_0^1(g_0)=[\quad], \, d_1^1(g_0)=[\quad], \, g \in \overline{W}_1 G,
s_0^0[\quad]=e_0,$$
Если $n \ge 1$, то:
$$d_0(g_{n-1},...,g_0)=(g_{n-2},...,g_0),$$
$$d_i(g_{n-1},...,g_0) = (d_{i-1} g_{n-1},...,d_0 g_{n-i} \cdot g_{n-i-1}, ..., g_0), \, 0 \le i \le n-1$$
$$d_n(g_{n-1}, ...,g_0) = (d_{n-1} g_{n-1}, ..., d_1 g_1),$$
$$s_0(g_{n-1},...,g_0) = (e_n, g_{n-1},..., g_0),$$
$$s_i(g_{n-1},...,g_0) = (s_{i-1} g_{n-1},...,s_0 g_{n-i},e_{n-i},g_{n-i-1},...,g_0).$$

Пусть $M$ - проконечный симплициальный $G$-модуль. Следуя Квиллену ~\cite{Qui1,Qui3,Qui5}, определим градуированную абелеву группу,
называемую гомологиями $G$ с коэффициентами в $M$, как
$$H_*^Q(G,M) := \pi_*(\mathbb{\widehat{Z}_{\mathfrak{C}}}[[W G]] \widehat{\otimes}_{\mathbb{\widehat{Z}_{\mathfrak{C}}}[[G]]} M),$$
где $\mathbb{\widehat{Z}_{\mathfrak{C}}}[[WG]]$ и $\mathbb{\widehat{Z}_{\mathfrak{C}}}[[G]]$ - свободные симплициальные абелевы группы,
полученные применением функтора группового кольца к $WG$ и $G$ соответственно, $\widehat{\otimes}$ - пополненное тензорное произведение,
см. \cite{ZR}.
Если мы рассмотрим в качестве группы $G$ постоянную группу, а в качестве $M$- постоянный
$\mathbb{\widehat{Z}}_{\mathfrak{C}}[[G]]$-модуль, то, как несложно увидеть, нормализованный комплекс
$\mathbb{\widehat{Z}_{\mathfrak{C}}}[[WG]]$ с дифференциалом $\partial_n = \sum_{i=0}^n (-1)^i d_i$ совпадает с неоднородной
бар-резольвентой \cite[6.2]{ZR}, и применением $\widehat{\otimes}_{\mathbb{\widehat{Z}_{\mathfrak{C}}}[[G]]} M$ получим  гомологии
$H_*(G,M)$ с коэффициентами в $M$ группы $G$. Ясно, что $\mathbb{\widehat{Z}_{\mathfrak{C}}}[[W
G]]\widehat{\otimes}_{\mathbb{\widehat{Z}_{\mathfrak{C}}}[[G]]} \mathbb{\widehat{Z}_{\mathfrak{C}}} =
\mathbb{\widehat{Z}_{\mathfrak{C}}}[[\overline{W} G]],$ где $\mathbb{\widehat{Z}_{\mathfrak{C}}}$ рассматривается как тривиальный
постоянный $G$-модуль, получаем, что $H_*(G,\mathbb{\widehat{Z}_{\mathfrak{C}}})\cong
\pi_*(\mathbb{\widehat{Z}_{\mathfrak{C}}}[[\overline{W} G]] ).$

Можно заменить проконечное кольцо $\mathbb{\widehat{Z}_{\mathfrak{C}}}$ на проконечное кольцо $\Lambda$ (например $\mathbb{F}_p$), см.
~\cite{ZR}. Тогда получаем $\Lambda[[W G]]$ - неоднородную бар-резольвенту с дифференциалами $\partial_n = \sum_{i=0}^n (-1)^i d_i,$ с
гомологиями $H_n(G,M) = Tor_n^{\Lambda[[G]]}(\Lambda,M)$.
В частности, в случае про-р-групп:
\begin{equation}
H_*(G,\mathbb{F}_p) = Tor_n^{\mathbb{F}_p[[G]]}(\mathbb{F}_p,\mathbb{F}_p) = \pi_*(\mathbb{F}_p[[W G]]
\widehat{\otimes}_{\mathbb{F}_p[[G]]} \mathbb{F}_p) = \pi_*(F_p[[\overline{W} G]]) = H_*^Q(G, \mathbb{F}_p)\label{8}
\end{equation}
В результате получаем:
\begin{predl}
Гомологии проконечной (про-р-группы) $H_*(G,\mathbb{\widehat{Z}}_\mathfrak{C})$ ($H_*(G,\mathbb{F}_p)$) совпадают с гомологиями  Квиллена
$H_*^Q(G,\mathbb{\widehat{Z}}_\mathfrak{C})$ ($H_*^Q(G,\mathbb{F}_p)$).
\end{predl}

\begin{opr}
Бисимплициальной про-$\mathfrak{C}$-группой называется совокупность

$\{G_{pq}, d_i^h, s_i^h, d_i^v, s_i^v,p,q \le 0\}$ про-$\mathfrak{C}$-групп $G_{pq}$, граничных непрерывных гомоморфизмов $d_i^h, d_i^v$
и непрерывных гомоморфизмов вырождения $s_i^h, s_i^v $, где $h$ - горизонтальная структура и $v$ - вертикальная структура,
удовлетворяющих следующим условиям:
$$d_j^h d_i^v= d_i^v d_j^h, \, s_j^h d_i^v = d_i^v s_j^h, \, s_j^v d_i^h = d_i^h s_j^v, \, s_j^v s_i^h = s_i^h s_j^v.$$
Каждая колонка и каждая строка бисимплициальной про-$\mathfrak{C}$-группы является симплициальной про-$\mathfrak{C}$-группой.
\end{opr}
Определим для заданной симплициальной про-$\mathfrak{C}$-группы $G$ симплициальную про-$\mathfrak{C}$-группу $EG$ следующим образом:
$$(EG)_n = \{ x \in G_{n+1}: d_0^1 d_0^2 ... d_0^{n+1}(x) = 1\}$$
с гомоморфизмами граней и вырождения, индуцированными из $G$, а также постоянную симплициальную про-$\mathfrak{C}$-группу $(C\pi_0G)_n =
\pi_0(G)$ с тождественными гомоморфизмами граней и вырождения.

Имеет место точная последовательность симплициальных про-$\mathfrak{C}$-групп
\begin{equation}
1 \to \Omega G \xrightarrow{i} EG \xrightarrow{\theta} G \xrightarrow{j} C\pi_0(G) \to 0 ,\label{5}
\end{equation}

где $\theta_n: (EG)_n \to G_n$ индуцирован гомоморфизмом $d_{n+1}^{n+1}: G_{n+1} \to G_n; \, j_n: G_n \to (C\pi_0 G)_n$ - это композит
$d_0^1 d_0^2... d_0^n$ с естественным гомоморфизмом $G_0 \to \pi_0(G)$ и $i: \Omega G \to EG$ - это гомоморфизм вложения ядра $\theta$.

Так как $\Im \theta = \Ker j$, а гомоморфизмы $j$ и $\theta$ являются расслоениями, то \ref{5} влечет следующие точные последовательности
гомотопических групп
\begin{eqnarray}
...\to \pi_n(\Omega G) \to \pi_n (EG)\to \pi_n (\Im \theta) \to \pi_{n-1}(\Omega G) \to... \label{16}\\
...\to \pi_n(\Im \theta) \to \pi_n(G) \to \pi_n(C \pi_0 G) \to \pi_{n-1}(\Im \theta) \to ...\label{6}
\end{eqnarray}

Мы утверждаем, что симплициальная про-$\mathfrak{C}$-группа $EG$ стягиваема справа со стягивающей гомотопией $h=0$ и $h_n =
s_{n+1}^{n+1}$. Действительно,
$\pi_0(EG):= \frac{(EG)_0}{d_1^1(\Ker d_0^2)} = \frac{\{x \in G_1 | d_0^1(x) = 1\}}{\{ d_1^1 (y) | y \in \Ker d_0^2\}}.$
$\pi_0(EG) = 0 \Leftrightarrow  (EG)_0 = \{ d_1^1 (y) | \text{ где } y \in \Ker d_0^2 \}$. Пусть $x \in G_1: d_0^1(x) = 1$, тогда
$s_1^1(x)$ - элемент в $G_2$ и $d_1^2 s_1^1(x) = x$. Надо проверить, что $d_0^2 s_1^1(x) = 1;$ действительно $d_0^2 s_1^1 (x) = s_0^0
d_0^1 (x) = 1 \Rightarrow \pi_0(EG)=0.$
Формулы $d_{n+1}^{n+1} h_n = 1$- это симплициальное тождество (10).
Надо еще проверить, что $0 = d_0^1 h_0(x) = h \varepsilon = 0$, но $d_0^1 s_1^1(x) = s_0 d_0^1(x) = 1$, аналогично проверяется и
последнее условие стягиваемости.

Из \ref{16},\ref{6} и стягиваемости $EG$  следует, что
\begin{eqnarray}
\pi_{n-1}(\Omega G) \cong \pi_n(G), n \ge 1 \label{7}
\end{eqnarray}

Определим диагональную симплициальную про-$\mathfrak{C}$-группу $\Delta G$:
$$(\Delta G)_n = G_{nn}, \quad d_j = d_j^h d_j^v,\quad s_j = s_j^h s_j^v.$$
Ясно, что если задана точная последовательность бисимплициальных про-$\mathfrak{C}$-групп
$0 \to G' \to G \to G'' \to 0$, то ей будет соответствовать длинная точная последовательность гомотопических групп: $$\to \pi_n(\Delta
G') \to \pi_n(\Delta G) \to \pi_n(\Delta G'') \to \pi_{n-1}(\Delta G') \to $$

Зафиксировав первый индекс в бисимплициальной про-$\mathfrak{C}$-группе $G$, получим вертикальную колонку, которая является
симплициальной про-$\mathfrak{C}$-группой. А с помощью \ref{5} точную последовательность бисимплициальных про-$\mathfrak{C}$-групп:
$$0 \to \Omega_v G_{**} \xrightarrow{i_v} E_vG_{**} \xrightarrow{\theta_v} G_{**} \xrightarrow{j_v} C_v \pi_0^v G_{**} \to 0.$$
$C_v \pi_0^v G$ - это бисимплициальная про-$\mathfrak{C}$-группа, у которой $(C_v \pi_0^v G)_{n,*} = \pi_0^v(G_{n,*}^v)$, горизонтальные
морфизмы  - это тождественные отображения группы $\pi_0 G_{n,*}$, а вертикальные индуцированы ${d_{n,0}^h}^i: G_{n,0} \to G_{n-1,0}.$

Тогда имеем точную последовательность диагональных про-$\mathfrak{C}$-групп
\begin{eqnarray}
0 \to \Delta(\Omega_v G_{*,*}) \to \Delta(E_v G_{*,*}) \to \Delta(G_{*,*}) \to \Delta ((C_v \pi_0^v G)_{*,*}) \to 0 \label{8}
\end{eqnarray}

$\Delta (C_v \pi_0^v G)$ - это симплициальная про-$\mathfrak{C}$-группа, у которой $\Delta (C_v \pi_0^v G)_n = \pi_0 (G_{n,*})$  и $d_i:
\Delta (C_v \pi_0^v G)_n \to \Delta (C_v \pi_0^v G)_{n-1} $ индуцирован ${d_{n,0}^h}^i: G_{n,0} \to G_{n-1,0}$

Нетрудно убедиться, что определенные таким образом гомоморфизмы граней и вырождений определены корректно и подчиняются необходимым
симплициальным тождествам, т.е. $\Delta C \pi_0 G$ действительно симплициальная про-$\mathfrak{C}$-группа. Ее группу гомологий комплекса
Мура в размерности $n$ будем обозначать $\pi_n^h \pi_0^v(G)$, т.е. $\pi_n^h \pi_0^v(G)= \pi_n(\Delta C_v \pi_0^v G)$

\begin{lemma}
$\Delta E_v G$ - стягиваемая симплициальная про-$\mathfrak{C}$-группа.
\end{lemma}
\textbf{Доказательство.}
Пользуясь перестановочностью $s^v, d^v$ и $s^h, d^h$  в определении бисимплициальной группы легко проверить, что гомоморфизмы $h=1$ и
$h_{nn}={s^v}_{n+1}^{n+1} {s^h}_{n+1}^{n+1}$ конструируют стягиващую справа гомотопию симпициальной про-$\mathfrak{C}$-группы $\Delta E_v
G.\boxtimes$

\begin{theorem}~\cite{Qui2}
Пусть $G$ - бисимплициальная про-$\mathfrak{C}$-группа, тогда имеются две спектральные последовательности:
$E_{nm}^2 = \pi_n^h \pi_m^v G \Rightarrow \pi_{n+m}(\Delta G),$
$E_{nm}^2 = \pi_n^v \pi_m^h G \Rightarrow \pi_{n+m}(\Delta G).$
\end{theorem}
\textbf{Доказательство.}
Последовательности \ref{8} соответствуют длинные точные последовательности \ref{16}, \ref{6}
$$...\to \pi_n(\Delta \Omega_v G) \to \pi_n(\Delta E_v G) \to \pi_n(\Im \Theta_v) \to \pi_{n-1}(\Delta \Omega_v G) \to...$$
 $$...\to \pi_n(\Im \Theta_v) \to \pi_n(\Delta G) \to \pi_n(\Delta C_v \pi_0^v G) \to \pi_{n-1}(\Im \Theta_v) \to ...$$
Получаем точную последовательность

\begin{equation}
\to \pi_{n-1} (\Delta \Omega_v G) \to \pi_n(\Delta G) \to \pi_n^h \pi_0^v G \to \pi_{n-2}(\Delta \Omega_v G) \to
\end{equation}

Обозначим $\Omega_v^m G:= \underbrace{\Omega_v \Omega_v ... \Omega_v}_{\mbox{$m$-раз}} G.$  Тогда из \ref{7}
получаем, что $\pi_0^v(\Omega_v^m G) = \pi_m^v(G)$ и, подставляя в \ref{8} вместо бисимплициальной про-$\mathfrak{C}$-группы $G$
бисимплициальную про-$\mathfrak{C}$-группу $\Delta \Omega_v^m G$, получаем точную пару
$$
\xymatrix{
{\pi_{n-1}(\Delta\Omega_v^{m+1}G)} \ar[rr] && {\pi_n(\Delta \Omega_v^m G)} \ar[ld] \\
& {\pi_n^h \pi_m^v G} \ar@{-->}[ul] &
}
$$
Для того, чтобы получить $r$, для которого $E_{nm}^r=\pi_n^h \pi_m^v(G)$, рассмотрим $d_r$ в виде композиции:
$$\pi_n^h \pi_m^v G \to \pi_{n-2}(\Delta \Omega_v^{m+1} G) \to \pi_{n-2}^h \pi_{m+1}^v(G)$$
откуда $r=2$ и $D_{m,n}^2 = \pi_n(\Delta \Omega_v^{m+1})$. Тогда стандартные аргументы о сходимости бирегулярной спектральной
последовательности, ~\cite[гл.~8]{Hu}, позволяют утверждать, что спектральная последовательность сходится к $\pi_{n+m}(\Delta G)$ с
фильтрацией
$F_n \pi_{n+m}(\Delta G) = \Im \{\pi_n(\Delta \Omega_v^m G)) \to \pi_{n+m}(\Delta G)\}.\boxtimes$

\section{Производные функторы.}
Одним из важнейших результатов агебраической топологии 70-х - это построение аксиоматической теории гомотопий \cite{Qui1}. Разработанная
Квилленом теория основывается на введении понятия модельной категории, т.е. категории с выделенными тремя классами морфизмов (расслоения,
корасслоения и слабые эквивалентности), играющими методологически ту же роль в теории гомотопий, что и давление, объем и температура в
молекулярной физике. Модельные категории с современной точки зрения приводят не только к унификации языка в теории гомотопий и
возможности адекватного сравнения различных теорий, но и, что особенно важно с нашей точки зрения, являются естественным вместилищем
неабелевой гомологической алгебры.

Для получения новых инвариантов в алгебраической геометрии Артин и Мазур \cite{AM} определили этальный гомотопический тип схемы, как
про-объект в гомотопической категории $\mathcal{H}$ категории симплициальных множеств $SS$. Затем Ректор \cite{Rec2} сделал определение
Артина-Мазура "жёстким", построив проконечное пополнение связного пространства, как про-объект в категории связных конечных
симплициальных множеств. Совершенно новый взгляд на ситуацию открылся благодаря работе Мореля \cite{M}, в которой были существенно
улучшены предыдущие результаты о р-пополнениях с помощью построения про-р-модельной структуры, где р - фиксированное простое число. В
\cite{M} была рассмотрена категория $SS^{profinite}$ - симплициальных проконечных пространств с модельной структурой, в которой слабые
эквивалентности - это отображения, индуцирующие изоморфизмы непрерывных $\mathbb{Z}/p$-когомологий. В частности, было показано, что
каждое проконечное пространство слабо эквивалентно в $SS^{profinite}$ проективному пределу конечных р-пространств с конечным набором
нетривиальных гомотопических групп, каждая из которых - конечная р-группа. Такая конструкция вместе с функтором  $SS \rightarrow
SS^{profinite}$ позволяет сделать "жёсткими" конструкции Артина-Мазура \cite{AM} для про-р-пополнений, а также обобщить работу Ректора на
случай неодносвязных пространств.
	В работе \cite{Q1}  Куик ввел несколько иную модельную структуру, в которой слабые эквивалентности - это отображения, индуцирующие
изоморфизмы на проконечных фундаментальных группах, а также в непрерывных когомологиях с конечными локальными коэффициентами, что
позволило определить не только $\mathbb{Z}/p$-пополнение для фиксированного простого числа p, как у Мореля, но и реализовать
про-$L$-модель в случае, когда $L$ содержит более чем одно простое число \cite{Q2}.

 Отметим похожий способ построения про-$L$-пополнений в работе \cite{Pri1}, а также элегантный обзор и сравнение пополнений по
 Бусфилду-Кану и $\mathbb{Z}/p$-пополнений в \cite{G}. Современное состояние  исследований про-гомотопических типов можно найти в работе
 Придхэма \cite{Pri1}.

Целью и задачей алгебраической топологии является сведение топологической проблемы к алгебраической, иными словами, построение
алгебраической модели. Для категории редуцированных множеств $SS_0$ Кан \cite{Kan,GJ} построил алгебраическую модель симплициальные
группы $SGr$, что на языке гомотопической алгебры означает Квиллен-эквивалентность \cite{GJ,Qui1} соответствующих гомотопических
категорий. Тогда, с точки зрения результатов Кана, обычная комбинаторная теория групп (задание групп образующими и соотношениями) ни что
иное, как двумерная теория гомотопий. Похожая ситуация имеет место и в проконечном случае, поскольку, как мы увидим в этой части работы,
существует та же пара сопряжённых функторов:

$$\xymatrix{ {SS_0^{profinite}} \ar@<1ex>[r]^{\widehat{Г}} & {SGr^{pro-\mathfrak{C}}} \ar@<1ex>[l]^{\overline{W}}},$$
обладающих рядом свойств, аналогичных дискретной ситуации. Это даёт основание определить место комбинаторной теории проконечных групп,
как двумерной проконечной теории гомотопий. Вернёмся к обсуждению данной аналогии уже в следующей главе путём развития теории
2-гомотопического типа, а именно, его алгебраической модели - проконечных скрещенных модулей.

\begin{predl}
Категория симплициальных про-$\mathfrak{C}$-групп является замкнутой симплициальной модельной категорией, у которой:

1) слабая эквивалентность $f: X \xrightarrow{\sim} Y$ - это гомоморфизм симплициальных про-$\mathfrak{C}$-групп, который индуцирует
изоморфизм $f_*: \pi_*(X) \xrightarrow{\sim} \pi_*(Y)$.

2) $f: X \twoheadrightarrow Y$ - расслоение, если $f$ - сюръективный гомоморфизм симплициальных про-$\mathfrak{C}$-групп.

3) $f: X \hookrightarrow Y$ - корасслоение, если $f$ имеет свойство подъема слева LLP, т.е. для любой коммутативной диаграммы
$$
\xymatrix{X \ar[r] \ar@{^(->}[d]^{f} & A \ar@{->>}[d]^{p}_{\sim} \\
Y \ar[r] \ar@{-->}[ur] & B},
$$
где $p$ - тривиальное расслоение (расслоение, которое также является и слабой эквивалентностью), существует пунктированная стрелка и
новая диаграмма тоже коммутативна.
\end{predl}

\textbf{Доказательство:}
Состоит в аккуратной проверке условий следующей теоремы:

\textbf{ Теорема} \cite[теор.4, гл. 2, \S 4]{Qui1}.
{\itshape Пусть $\mathcal{A}$ - категория, замкнутая относительно конечных прямых пределов и имеющая достаточно много проективных
объектов. Пусть $s\mathcal{A}$ - симплициальная категория симплициальных объектов над $\mathcal{A}$. Определим отображение $f$ в
$s\mathcal{A}$ как расслоение (соотв. слабая эквивалентность), если $\underline{Hom}(P,f)$ (под $P$ понимается постоянный симплициальный
объект)  является расслоением (соотв. слабой эквивалентностью) в категории симплициальных множеств для каждого проективного объекта $P$ в
$\mathcal{A}$, и корасслоение, если $f$ имеет LLP относительно тривиальных расслоений. Тогда $s\mathcal{A}$ - замкнутая симплициальная
модельная категория, если $\mathcal{A}$ удовлетворяет одному из следующих дополнительных условий:

(*) Каждый объект в $s\mathcal{A}$ является фибрантным.

(**) $\mathcal{A}$ замкнута относительно индуктивных пределов и имеет множество малых проективных образующих.}

 Если $\mathcal{A}$ - категория про-$\mathfrak{C}$-групп, то условие *) выполнено, т.к. каждый объект $G$ в $\mathcal{A}$ является
 фактором когруппового объекта $C$ \cite[Предл.1, гл.2, 4.2]{Qui1}. В качестве $C$ следует взять свободную про-$\mathfrak{C}$-группу
 $F_{\mathfrak{C}}(X)$ на пространстве образующих $G$.

 В случае категории симплициальных про-$\mathfrak{C}$-групп расслоения и слабые эквивалентности из теоремы  \cite[теор.4, гл. 2, \S
 4]{Qui1} совпадают с нужными. Действительно, пусть $\mathcal{A}$ - категория про-$\mathfrak{C}$-групп и $f: A \to B$ - некоторый морфизм
 в $s\mathcal{A}$, тогда для любого $C \in \mathcal{A}$ определяется
$\underline{Hom}(C,A)\in SS,$ полагая $\underline{Hom}(C,A)_n = Mor_\mathcal{A}(C, A_n)$, при этом морфизмы граней и вырождения
получаются из соответствующих морфизмов граней и вырождения в $s\mathcal{A}$.
Если в качестве $C$ взять свободную про-$\mathfrak{C}$-группу над проконечным пространством $X$, то
$Mor_\mathcal{A}(F_\mathfrak{C}(X),A_n) \cong \prod_X A_n$
имеет $(F_\mathfrak{C}(X)$ - когрупповой объект в $A$) структуру группы, как прямое произведение про-$\mathfrak{C}$-групп $A_n$.
Если вспомнить \cite[гл.2, 3.7-3.13]{Qui1}, что сюръективные гомоморфизмы симплициальных групп суть расслоения, а гомотопические группы
(как гомологии комплекса Мура) совпадают с гомотопическими группами базовых симплициальных множеств, то получаем нужное совпадение
модельных структур. Для произвольной проективной про-$\mathfrak{C}$-группы доказательство проводится с помощью ретракт-аргумента, см.
\cite[Предл.1, гл. ~2, стр. ~4.2]{Qui1}$\, \boxtimes$.

Следуя \cite{Kan,GJ}, см. также \cite{Q2}, построим функтор $\widehat{Г}^\mathfrak{C}$ из категории редуцированных симплициальных
проконечных пространств в категорию симплициальных про-$\mathfrak{C}$-групп. Во-первых, согласно \cite{Qui4,M,Q2} заметим, что любая
симплициальная про-$\mathfrak{C}$-группа и любое симплициальное проконечное пространство есть проективный предел конечных симплициальных
групп и конечных симплициальных пространств соответственно. Если $K$ - редуцированное конечное симплициальное множество, тогда определим
$\widehat{Г}_n^\mathfrak{C}(K),$ полагая $\widehat{Г}_n^\mathfrak{C}(K):=F_\mathfrak{C}(K_{n+1}-s_0(K_n)).$
Для произвольного проконечного редуцированного симплициального пространства $K$: $$\widehat{Г}_n^\mathfrak{C}(K):=\varprojlim_{\la \in
\Lambda} F_{\mathfrak{C}}(({K_{\la}})_{n+1}-s_0(K_{\la})_n) \cong F_{\mathfrak{C}}(K_{n+1}-s_0(K_n))$$

Пусть $x \in K_{n+1},$ обозначим $\tau(x)$ - класс элемента $x$ в $\widehat{Г}_n^\mathfrak{C}(K)$, тогда операторы граней и вырождения в
$\widehat{Г}_n^\mathfrak{C}(K)$ задются следующим образом:
$$\tau(d_0x) d_0 \tau(x) = \tau(d_1 x)$$
$$d_i \tau(x) = \tau(d_{i+1} x), \, i > 0$$
$$s_i \tau(x) = \tau (s_{i+1} x), \, i \ge 0$$
Определим теперь $d_i$ и $s_i$ путем единственного расширения до гомоморфизмов про-$\mathfrak{C}$-групп. Немедленно проверяется, что
$\widehat{Г}^\mathfrak{C}(K)$ является комплексом групп, и $\tau: K \to \widehat{Г}^\mathfrak{C}(K)$ - так называемая скручивающая
функция \cite{May}.

Заметим, следуя \cite{Kan,GJ}, что функторы $\widehat{Г}^\mathfrak{C}$ и $\overline{W}$ сопряжены. При этом, изоморфизм
 $$Mor_{SGr^{prof}}(\widehat{Г}^\mathfrak{C}(K), A) \cong Mor_{SS^{prof}}(K, \overline{W}A)$$ строится следующим образом. Пусть задан
 гомоморфизм проконечных симплициальных групп $f: \widehat{Г}^\mathfrak{C}(K) \to A$ и пусть $f(\tau(x)) = \sigma$. Определим
$$\Psi(f): K \to \overline{W}A$$ по правилу $$\Psi(f)(x) = (f(\tau(x)), \, d_0 f(\tau(x)),...,d_0^n f(\tau(x))$$
Обратно, если задан некоторый морфизм симплициальных проконечных пространств $$g: K \to \overline{W}A,$$
$$g(x) = [\tau_{n-1}, ..., \tau_0],$$
тогда определим
$$\Phi(g): \widehat{Г}^\mathfrak{C}(K) \to A$$
по правилу $$\Phi(g)(\tau(x)) = \tau_{n-1}$$
Прямая проверка показывает, что $\Psi \Phi(g) = g$ и $\Phi \Psi(f) = f,$ а также, что $\Phi(g)$ и $\Psi(f)$ являются симплициальными
гомоморфизмами.

\textbf{ Пусть далее до конца этой части работы $\mathfrak{C}$ -класс всех конечных групп, порядки которых суть произведения простых
чисел $l$ из заданного набора простых чисел $L$, тогда единица сопряжения, т.е. гомоморфизм симплициальных групп $\Phi \Psi(id_G):
\widehat{Г}^\mathfrak{C} \overline{W} G \to G,$ где $G$ - симплициальная прo-$\mathfrak{C}$-группа, является слабой гомотопической
эквивалентностью.}

Разложим симплициальную про-$\mathfrak{C}$-группу в проективный предел конечных симплициальных $\mathfrak{C}$-групп $G = \varprojlim
G_{\la}$ и воспользуемся эквивалентностью Кана \cite{Kan,GJ} для  конечной симплициальной группы $G_{\la}$. Имеем $\overline{W}G_{\la}$ -
конечные симплициальные пространства, у которых $\pi_{n+1} \overline{W} G_{\la} \cong \pi_n G_{\la},$ а $\overline{W}G \cong \varprojlim
\overline{W} G_{\la}.$ Функтор Кана $Г$ (дискретный) ставит в соответствие конечному симплициальному множеству $\overline{W}Г_{\la}$
дискретную симплициальную группу
   $Г\overline{W}G_{\la}$, причём имеем:
   $\pi_n(Г\overline{W}G_{\la}) \cong \pi_{n+1}(\overline{W}G_{\la}) \cong \pi_n(G_{\la}),$
но $Г\overline{W}G_{\la}$ - конечнопорожденные свободные группы (дискретные), а следовательно \cite{Qui3}, \cite{Rec2} они являются
"хорошими" и "$L$-хорошими", и $\pi_n(\widehat{Г}^\mathfrak{C} \overline{W}G_{\la}) \cong
\widehat{\pi_n(Г\overline{W}G_{\la})^\mathfrak{C}}.$
Ввиду точности проективных пределов в категории проконечных пространств: $\pi_n(\widehat{Г}^\mathfrak{C} \overline{W}G_{\la}) \cong
\varprojlim \pi_n(\widehat{Г} \overline{W} G_{\la}) \cong \varprojlim \widehat{\pi_n(Г \overline{W} G_{\la})}^\mathfrak{C} \cong
\varprojlim \pi_n(G_{\la}) \cong \pi_n(G).$ Т.к. $\widehat{Г}^\mathfrak{C} \overline{W}G \cong \varprojlim \widehat{Г}^\mathfrak{C}
\overline{W} G_{\la}$, то переходом к проективному пределу, получаем:
\begin{predl}
Пусть $\mathfrak{C}$- класс конечных групп, порядки которых - произведения степеней простых чисел из заданного набора $L$; $G$ -
симплициальная про-$\mathfrak{C}$-группа, тогда естественный гомоморфизм $$\Phi\Psi(id_G): \widehat{Г}^\mathfrak{C} \overline{W} G \to
G$$ является слабой эквивалентностью симплициальных про-$\mathfrak{C}$-групп.
\end{predl}

Представим симплициальное редуцированное проконечное пространство в виде проективного предела конечных симплициальных пространств
$\overline{W}G \cong \varprojlim K_\lambda$, $(K_\lambda)_n = \bigcup s_i((K_\lambda)_{n-1})\biguplus (Y_\lambda)_n$, Положим $K_n=
\bigcup s_i(K_n)\biguplus \varprojlim Y_\lambda$, тогда свободная симплициальная про-$\mathfrak{C}$-группа
$\widehat{Г}^\mathfrak{C}\overline{W}G$ имеет CW-базис и кофибрантна, т.к. морфизм $* \hookrightarrow
\widehat{Г}^\mathfrak{C}\overline{W}G $ является почти свободным гомоморфизмом, см. \cite[лемма. ~6.1]{GJ}. Следовательно, в нашем
распоряжении имеются функториальные кофибрантные замены в категории симплициальных про-$\mathfrak{C}$-групп. Нетрудно заметить, что
симплициальные резольвенты, построенные в первой части работы, также кофибрантны.

Теперь заметим, что свободная симплициальная про-$\mathfrak{C}$-резольвента может быть построена функториально, как
$\widehat{Г}^\mathfrak{C}\overline{W}G$, где $\widehat{Г}^\mathfrak{C}$ - про-$\mathfrak{C}$-функтор Кана, $\overline{W}$ - функтор
классифицирующего пространства, а $G$ - постоянная симплициальная про-$\mathfrak{C}$-группа, у которой $G_n = G$, а морфизмы граней и
вырождения - тождественные гомоморфизмы $G$. Действительно, $\widehat{Г}^\mathfrak{C} \overline{W} G$ - свободная симплициальная группа,
c другой стороны предыдущие рассуждения показывают, что $\widehat{Г}^\mathfrak{C} \overline{W} G$ слабо эквивалентна $G$, но $\pi_n(G) =
0$ при $n \ge 1$, а $\pi_0(G)=G.$

Также сформулируем следствие, которое нам понадобится в последней части работы.
\begin{sled}
Пусть $G$ - конечная $\mathfrak{C}$-группа, тогда существует конечнопорожденная свободная симплициальная про-$\mathfrak{C}$-резольвента
$G$.
\end{sled}

Для модельных категорий естественным образом строятся их гомотопические категории \cite{GJ,Qui1}, которые по построению представляют
собой категорию $Ho(\mathcal{C})$, у которой объекты те же, что и в исходной категории $\mathcal{C}$, а морфизмы - гомотопические классы
$Hom_{\pi C_{cf}}$ $$Hom_{H_0(\mathcal{C})}(X,Y) = Hom_{\pi C_{cf}}(R'QX,R'QY) = \pi(RQX,RQY),$$
где $QX,QY$ - кофибрантные замены объектов $X$ и $Y \in \mathcal{C}$, получаемые путем разложения отображений $\varnothing \to X$ в
композицию
$\xymatrix{\varnothing \ar[r] & QX  \ar@{->>}[r]^{p_x}& X}$, где $p_x$ - тривиальное расслоение; $RX$ получается из разложения
отображения $X \to *$ в композицию \xymatrix{X \ar@{^(->}[r]^{i_x} & RX \ar@{->>}[r] & {*}}, где $i_x$ - тривиальное корасслоение.

 В случае модельных категорий, $Ho(\mathcal{C})$ имеет привычную интерпретацию \cite{GZ}, как локализация $\mathcal{C}$ относительно
 морфизмов, которые являются слабыми эквивалентностями. Более точно:
\begin{opr}\cite{GZ}
Пусть $\mathcal{C}$ - категория, $W \subseteq Mor(\mathcal{C})$ - класс морфизмов. Функтор $F: \mathcal{C} \to \mathcal{D}$ называется
локализацией $\mathcal{C}$ относительно $W$, если:

(i) $F(f)$ является изоморфизмом для каждого $f \in W,$

(ii) если $G: \mathcal{C} \to \mathcal{D}'$ - функтор, переводящий элементы из $W$ в изоморфизмы, тогда существует единственный функтор
$G': \mathcal{D} \to \mathcal{D}', G'F = G.$
\end{opr}
\begin{theorem}\cite[гл.1, 1.13]{Qui2},\cite{GJ}
Пусть $\mathcal{C}$ - модельная категория и $W$ - класс слабых эквивалентностей. Тогда функтор $\gamma: \mathcal{C} \to Ho(\mathcal{C})$
является локализацией $\mathcal{C}$ относительно $W$.
\end{theorem}

Гомологическую алгебру для модельных категорий удобно определить, как левые или правые производные функторы $LF,RF: Ho(\mathcal{C}) \to
\mathcal{D}$ для некоторого функтора $F: \mathcal{C} \to \mathcal{D}$.
\begin{opr}\cite[гл.1, \S 4]{Qui2}
Предположим, что $\mathcal{C}$ - модельная категория и что $F: \mathcal{C} \to \mathcal{D}$ - функтор. Рассмотрим пары $(G,s)$, состоящие
из функтора $G: Ho(\mathcal{C}) \to \mathcal{D}$ и естественного преобразования $s: G \gamma \to F$. Тогда левым производным функтором
для $F$ называется пара $(LF,t)$ универсальная слева в следующем смысле, что если $(G,s)$ - какая-то другая пара, то существует и
единственное естественное преобразование $s': G \to LF$ такое, что композит естественных преобразований $G \gamma \xrightarrow{s' \circ
\gamma} (LF) \gamma \xrightarrow{t} F$ является естественным преобразованием $s$.
(Аналогичным образом определяется правый производный функтор \cite{Qui2}.)
\end{opr}

Следующее утверждение-определение и его следствие с избытком покрывают необходимые нам основания для построения неабелевой гомологической
алгебры в категории про-$\mathfrak{C}$-групп:

\begin{opr}\cite[7.1]{GJ}
Пусть $\mathcal{C}$ - симплициальная модельная категория и $\mathcal{A}$ - некоторая категория. Предположим $F: \mathcal{C} \to
\mathcal{A}$ - функтор, который переводит слабые эквивалентности между кофибрантными объектами в изоморфизмы. Определим тотальный левый
производный функтор
$$ LF: Ho(\mathcal{C}) \to \mathcal{A},$$
полагая $LF(X)=F(Y)$, где $Y \rightarrow X$ - тривиальное расслоение и $Y$ кофибрантен.
\end{opr}

\begin{sled}
Пусть $X \in \mathcal{C}$. Если $X$ кофибрантен, то $LF(X) \cong F(X)$. Если $Y \rightarrow X$ слабая эквивалентность и $Y$ кофибрантен,
то $LF(X) \cong F(Y).$
\end{sled}

Сконцентрируем внимание на функторах $F: sGr^{pro-\mathfrak{C}} \to Ab^{pro-\mathfrak{C}}$, которые могут быть представлены в виде
композиции $F=AT$, где
$$ T: sGr^{pro-\mathfrak{C}}  \to sGr^{pro-\mathfrak{C}},$$
при этом $T$ сохраняет симплициальные гомотопические эквивалентности и
$$ A: sGr^{pro-\mathfrak{C}} \to Ab^{pro-\mathfrak{C}} $$
функтор гомотопических групп симплициальной про-$\mathfrak{C}$-группы
$$ A(G)=\prod_{i \geq 0} \pi_i(G).$$

Функторы такого вида удовлетворяют условиям определения 11, т.к. $T$ переводит слабые эквивалентности между кофибрантными объектами в
слабые эквивалентности. Заметим, каждый объект в $sGr^{pro-\mathfrak{C}}$ фибрантен. По теореме Уайтхеда ~\cite[теорема 1.10]{GJ}, если
$f$ - слабая гомотопическая эквивалентность между кофибрантными объектами, то $f$ - гомотопическая эквивалентность.  Для построения
гомотопии в доказательстве теоремы Уайтхеда можно выбрать канонический цилиндр-объект вида $X \otimes I$ ( т.к. у нас симплициальная
модельная категория); поэтому можно считать, что $f$ - симплициальная гомотопическая эквивалентность. Тогда, т.к. $T$ сохраняет
симплициальные гомотопические эквивалентности, то $T(f)$ - симплициальная гомотопическая эквивалентность. Откуда по ~\cite[предложение
1.14]{GJ} $T(f)$ - слабая эквивалентность. Очевидно, что $A$ переводит слабые эквивалентности в изоморфизмы.

Множество осмысленных примеров возникает из функторов $T: Gr^{pro-\mathfrak{C}} \to Gr^{pro-\mathfrak{C}}$, обладающих свойством
$T(id)=id$, т.к. последнее условие обеспечивает сохранение симплициальных гомотопических эквивалентностей функтором $sT:
sGr^{pro-\mathfrak{C}} \to sGr^{pro-\mathfrak{C}}$, полученным покомпонентным применением $T$ к симплициальной группе.

Рассмотрим, следуя линии \cite[A.13]{Mikh}, примеры функторов, представимых в виде $A \circ sT$. Пусть
$$T:  Gr^{pro-\mathfrak{C}} \to Gr^{pro-\mathfrak{C}}, \, T(1) = 1,$$
тогда из наших предыдущих рассуждений следует, что $\pi_n(T(X))$ не зависит от выбора свободной симплициальной
про-$\mathfrak{C}$-резольвенты $X$ про-$\mathfrak{C}$-группы $G$ (это следует также и из [теор. ~1.5]), и определим левые производные
функторы функтора $T$, как
$$\mathfrak{L}_i(T): G \to \pi_i(T(X)), \, i \ge 0.$$
Так как $\pi_i(T(X))$ - абелевы группы при $i \ge 1$ \cite{May}, то получаем функторы $$\mathfrak{L}_i(T): Gr^{pro-\mathfrak{C}} \to
AbGr^{pro-\mathfrak{C}}, \, i \ge 1$$ $$\mathfrak{L}_0(T):  Gr^{pro-\mathfrak{C}} \to Gr^{pro-\mathfrak{C}}, \, i=0$$

\begin{prim}[Производные функтора mod-p-абелианизации]
Пусть $G$ - про-р-группа, обозначим $\Phi(G)$ ее подгруппу Фраттини. Тогда
$$Ab_p: Gr^{pro-p} \to AbGr^{pro-p}$$
- функтор из категории про-р-групп в категорию абелевых про-р-групп экспоненты р, который ставит про-р-группе ее фактор-группу по
подгруппе Фраттини, которая \cite{Se1} для заданной про-р-группы $G$ совпадает с $\overline{G^p[G,G]}$.

Пусть $X$ - свободная симплициальная резольвента про-р-группы $G$. Мы утверждаем, что имеет место следующее:
\begin{predl}
$\mathfrak{L}_nAb_p(G):=\pi_n(X/\Phi(X)) \cong H_{n+1}(G), \text{ где } H_n(G)=H_n(G,\mathbb{F}_p) \text{ в смысле \cite{ZR}}.$
\end{predl}
\end{prim}
\textbf{Доказательство.}
Рассмотрим бисимплициальную про-р-группу:

$Y_{i,j} = \mathbb{F}_p[[\overline{W}X_i]_j] = \mathbb{F}_p[[\underbrace{X_i \times ... \times X_i}_{\mbox{$j$-раз}}]]$, в которой
вертикальная симплициальная структура (на пополненном групповом кольце) возникает из симплициальной структуры на $X_i$, а горизонтальная
исходит из симплициальной структуры на $\mathbb{F}_p[[\overline{W}X]]$, так как $\overline{W} X$ - симплициальное проконечное
пространство.

1) Если зафиксировать индекс $j$, то мы получим симплициальную группу вида $\mathbb{F}_p[[\underbrace{X \times ... \times
X}_{\mbox{$j$-раз}}]]$, $i$ изменяется. Так как $X_i$ - свободная симплициальная резольвента, то:
$$\pi_0(\mathbb{F}_p[[\underbrace{X \times ... \times X}_{\mbox{$j$-раз}}]]) \cong \mathbb{F}_p[[\overline{W}_j(G)]], \, i \ge 0, \quad
\pi_n(\mathbb{F}_p[[\underbrace{X \times ... \times X}_{\mbox{$j$-раз}}]]) = 0, \, n \ge 1$$

Это легко следует из описания ядер гомоморфизмов групповых колец, индуцированных гомоморфизмами соответствующих групп ~\cite{GSh,KH}.
Замечая, что гомологии комплекса $\mathbb{F}_p[[\overline{W}_j (G)]]$ совпадают с гомологиями про-р-группы $G$, получаем
$$\pi_n^v \pi_0^h (Y_{ij}) = H_n(G).$$

Будем использовать спектральную последовательность Квиллена со сдвигом вверх на единицу. А именно
$$E_{n,m}^2 := \pi_n^v \pi_{m-1}^h (Y_{ij}), D_{n,m}^2 = \pi_n^v (\Delta \Omega_v^m G) $$
Для этой бирегулярной спектральной последовательности имеем
$E_{n,m}^2 = \begin{cases}
H_n(G),&\text{если $m=1$},\\
0,&\text{если $m \ne 1$.}
\end{cases}
$
С другой стороны, в обозначениях ~\cite[гл.~8]{Hu} $E_{n,1}^2\cong E_{n,1}^{\infty}=\frac{H_{n,1}}{H_{n-1,2}},$ где $H_{n,1} =
D_{n+2,-1}^3,$
а $H_{n-1,2} = D_{n+1,-1}^4.$ Мы утверждаем, что $H_{n-1,2} = 0$, действительно, в силу ~\cite[гл.~8, (P$\partial 8$)]{Hu}, имеем точную
последовательность
$0 \to H_{n-2,3} \to H_{n-1,2} \to E_{n-1,2}^{\infty} \to 0,$ но $E_{n-1,2}^{\infty} = 0$, следовательно $H_{n-2,3} \cong H_{n-1,2},$
аналогично для $H_{n-2,3} \quad$
$0 \to H_{n-3,4} \to H_{n-2,3} \to E_{n-2,3}^{\infty} \to 0,$ но $E_{n-2,3}^{\infty} = 0 \Rightarrow H_{n-3,4} = H_{n-2,3}.$
Итерируя этот спуск, получим на $(n-1)$-ом шаге, что $H_{0,n+1} \cong H_{-1,n+2} = 0$, откуда $H_{n-1,2} = 0.$

Собирая все вычисления вместе, получаем $E_{n,1}^2 \cong E_{n,1}^{\infty} \cong \frac{H_{n,1}}{H_{n-1,2}} \cong H_{n,1} \cong
D^3_{n+2,-1}.$

Для вычисления $D_{n+2,-1}^3$ рассмотрим последовательность точной пары
$$0 = E_{n+2,0}^2 \xrightarrow{k_2} D_{n+1,0}^2 \xrightarrow{i^2} D_{n+2,-1}^2 \xrightarrow{j^2} E_{n+1,0}^2, n \ge 0.$$
Т.к. $E^2_{n+2,0} = E_{n+1,0}^2 = 0$ получаем, что $i^2$ - изоморфизм, следовательно $D_{n+2,-1}^3 = i^2(D_{n+1,0}^2) \cong
D_{n+1,0}^2.$
Остается заметить, что $D_{n+1,0}^2 = \pi_{n+1}(\Delta \Omega^0 Y) = \pi_{n+1}(\Delta Y).$

2) Если зафиксировать  индекс $i$, то получим вертикальную симплициальную про-р-группу, которая суть $\mathbb{F}_p[[\overline{W} X_i]]$.
Заметим, что ее гомологии - это в точности гомологии свободной про-р-группы $X_i$ с коэффициентами в $\mathbb{F}_p,$ которые, как
известно ~\cite{ZR}, равны
$H_0(X_i) = \mathbb{F}_p, H_1(X_i) = X_i/\Phi(X_i)$ и $H_n(X_i) = 0$ для $i>1.$
Получаем, подставляя в спектральную последовательность Квиллена
$$E_{n,m}^2 = \pi_n^h \pi_m^v (Y_{ij}) = \begin{cases}
\pi_n(X/\Phi(X)),&\text{m=1},\\
\mathbb{F}_p,&\text{m=0,n=0},\\
0, &\text{если $m \ne 1$ и $(m=0,n=0)$}.
\end{cases}
$$
Для этой спектральной последовательности аналогично первому случаю
$$\pi_n(X/\Phi(X)) = E_{n,1}^2 = E_{n,1}^{\infty} \cong D_{n+1,0}^2 \cong \pi_{n+1}(\Delta \Omega_h^1 Y) \cong \pi_{n+2}(\Delta Y) $$
Соединяя случай 1) и случай 2), получаем
$\mathcal{L}_n Ab_p (G)=\pi_n(X/\Phi(X)) \cong H_{n+1}(G).\, \boxtimes$

\begin{prim}[Функтор Фраттини]
Пусть $G$ -  про-р-группа, $$T: Gr^{pro-p} \to Gr^{pro-p}$$  - функтор, заданный формулой $$T: G \to \overline{G^p[G,G]}$$

Если задано копредставление $1 \to R \to F \to G \to 1,$ то без ограничения общности можно считать, что $R \subset F^p[F,F].$ Этого
всегда можно добиться, исключая из соотношений образующие $F$ и уменьшая, если надо, ранг $F$, т.е. используя минимальное копредставление
про-р-группы $G$ $(rank F = dim_{\mathbb{F}_p} H^1(G))$, количество соотношений будет равно $dim_{\mathbb{F}_p}H_2(G)$, см.\cite{ZR}.

Покажем, что $\mathcal{L}_i \Phi(G) = \begin{cases}
\frac{\overline{F^p[F,F]}}{\overline{R^p[R,F]}},& i=0,\\
H_{i+2}(G),& i>0.
\end{cases}$

Рассмотрим короткую точную последовательность симплициальных про-р-групп $$1 \to \Phi(X) \to X \to X/\Phi(X) \to 1.$$
Т.к. $\pi_n(X)=0$ при $n \ge 1$, сразу получаем, что $\pi_n(\Phi(X)) = \pi_{n+1} (X/\Phi(X)), \, n \ge 1$, следовательно
$\mathcal{L}_i \Phi(X)  = H_{i+2}(G), \, i > 0.$

Вычислим значение $\mathcal{L}_0 \Phi(G)= \pi_0(\Phi(X))$.
В малых размерностях
$0 \to \pi_1(X/\Phi(X)) \to \pi_0(\Phi(X)) \xrightarrow{\varepsilon} \pi_0(X) \to \pi_0(X/\Phi(X)) \to 0.$
Из ранее доказанного $\pi_1(X/\Phi(X)) \cong H_2(G) \cong R/\overline{R^p[R,F]}$.
Из анализа связывающего гомоморфизма комплексов Мура симплициальных про-$\mathfrak{C}$-групп, полученных из некоторого копредставления
$G,$ любой элемент $(N(X/\Phi(X)))_1$ является образом некоторого элемента из $(NX)_1$ (доказательство Теоремы 1.4 в \cite{Ina}), но
$d_1^1(NX)_1$ совпадает с $R$. С другой стороны, т.к. $R\subset \Phi(X_0),$ ясно, что гомоморфизм $\pi_1(X/\Phi(X)) \to \pi_0(\Phi(X))$
совпадает с гомоморфизмом $\frac{R}{\overline{R^p[R,F]}} \to \pi_0(\Phi(X_0));$
при этом, т.к. нам известно, что ядро такого гомоморфизма нулевое, имеем
$1 \to \frac{R}{\overline{R^p[R,F]}} \to \pi_0(\Phi(X)) \xrightarrow{\varepsilon} G \to G/\overline{G^p[G,G]} \to 1$, откуда $1 \to
\frac{R}{\overline{R^p[R,F]}} \to \pi_0(\Phi(X)) \to \overline{G^p[G,G]} \to 1.$
С другой стороны $\overline{G^p[G,G]} \cong \frac{\overline{F^p[F,F]}}{R},$ т.к. $G/(\overline{F^p[F,F]}/R) \cong
(F/R)/(\overline{F^p[F,F]}/R) \cong F/\overline{F^p[F,F]} \cong G/\overline{G^p[G,G]}$,
тогда
$1 \to \frac{R}{\overline{R^p[R,F]}} \to \frac{\overline{F^p[F,F]}}{\\Im d_1^1(\Ker d_0^1)} \to \frac{\overline{F^p[F,F]}}{R} \to 1
\Rightarrow \Im d_1^1(\Ker d_0^1) = \overline{R^p[R,F]}.$

\end{prim}

\begin{prim}[Проконечный функтор Ли]

Пусть $\mathfrak{C}$-класс всех конечных групп, $T:Gr^{profinite} \to AbGr^{profinite}$  - функтор, заданный формулой $G \to
G/\overline{[G,G]}.$ Вычисления с использованием спектральной последовательности Квиллена, аналогичные примеру 3, дают
$$\mathcal{L}_i(Ab) \cong H_i(G,\mathbb{\widehat{Z}})$$

Несколько сложнее ситуация в случае, когда функтор $T: Gr^{profinite}\to Gr^{profinite}$ определяется правилом $$T: G \to
\overline{[G,G]}.$$

Имеем расслоение симплициальных про-$\mathfrak{C}$-групп
$1 \to \overline{[X,X]} \to X \to X/\overline{[X,X]} \to 0$

Соответствующая длинная точная последовательность гомотопических групп влечёт $\mathcal{L}_i(Lie G) \cong \mathcal{L}_{i+1}(Ab), i \ge 1,
\text{ т.е. } \mathcal{L}_i(Lie G) = H_{i+2}(G,\widehat{\mathbb{Z}}), i \ge 1.$

Остается посчитать $\mathcal{L}_0(Lie).$ Имеем точную последовательность
$$
0 \to \pi_1(X/\overline{[X,X]}) \to  \pi_0(\overline{[X,X]}) \to \pi_0(X) \to  \pi_0(X/\overline{[X,X]}) \to 1
$$

\end{prim}

Откуда
$
\xymatrix{
0 \ar[r]  & \pi_1(X/\overline{[X,X]}) \ar[r] & \pi_0(\overline{[X,X]}) \ar[r]^{\varepsilon} & \overline{[G,G]} \ar[r] & 0\\
}.$
По формуле Хопфа, см. ~\cite[~8.2.3]{P1} $\pi_1(X/\overline{[X,X]}) \cong H_2(G, \widehat{\mathbb{Z}}) =
\frac{R\cap\overline{[F,F]}}{\overline{[F,R]}}.$
Обозначим $F=X_0$, тогда если
$\varepsilon$ профакторизовать по $\overline{[F,R]}$, получим $0 \to \frac{R\cap\overline{[F,F]}}{\overline{[F,R]}} \to
\frac{\overline{[F,F]}}{\overline{[F,R]}} \to \overline{[G,G]} \to 0$, откуда  $\pi_0(Lie G) \cong
\frac{\overline{[F,F]}}{\overline{[F,R]}}
$

\section{От симплициальных резольвент к скрещенным и цепным комплексам (2-мерная проконечная гомотопия).}

Как было показано в предыдущей главе, свободная симплициальная про-$\mathfrak{C}$-резольвента $X$ занимает высшую позицию в иерархии
объектов, характеризующих гомотопические свойства про-$\mathfrak{C}$-группы $G$. Из первой части работы следует, что задание образующими
и соотношениями следует рассматривать, как некоторую 1-срезанную свободную симплициальную про-$\mathfrak{C}$-резольвенту.

  Напомним, что имеется пара сопряжённых функторов:
  $\xymatrix{sGr^{pro-\mathfrak{C}}\ar@<1ex>[r]^{tr^n} & sGr^{pro-\mathfrak{C}}_{tr^n}\ar@<1ex>[l]^{cosk^n},}$
  где функтор $cosk^n$- является правым сопряжённым к $tr^n$. Известно \cite [1.3]{CC} $\pi_n(cosk^1sk^1(X))=0$ при $n\geq~2$, тем самым,
  с гомотопической точки зрения, комбинаторная теория групп суть вопрос об изучении 2-типов или, как объяснено в \cite{CC,P1,P2}, их
  алгебраической модели - скрещенных модулей.

  С каждым копредставлением, следуя \cite{BH,P1,P2}, можно связать скрещенный модуль копредставления, который можно использовать для
  построения начальной части свободной $\mathbb{Z}_\mathfrak{C}[[G]]$-резольвенты, используя последовательность Кроуэла-Линдона. В
  фундаментальной работе о классификации теоретико-групповых типов дискретных симплициальных групп \cite{CC} Карраско и Цегарра для
  каждой симплициальной группы $H$ ввели комплекс $(C_n, \partial_n = \overline{d_n^n}),$ где $С_n(H) = \frac{N H_n}{(N H_n \cap D_n)
  d_{n+1} (NH_{n+1} \cap D_{n+1})},$
   построенный из комплекса Мура($N H_n, d_n^n$), $D_n$ - подгруппа в $H_n$, порожденная вырожденными элементами.
  В \cite{MP} доказано, что комплекс $\mathfrak{C}(H)$, являющийся комплексом скрещенных модулей, а в случае, когда $H$ - свободная
  симплициальная резольвента, будет свободной скрещенной резольвентой группы. Мы не ставим себе задачу перенести эти результаты в
  про-$\mathfrak{C}$-случай в полной общности, однако покажем, что скрещенный модуль, возникающий из свободной симплициальной
  резольвенты, согласно идеологии \cite{CC}, совпадает со скрещенным модулем про-$\mathfrak{C}$-копредставления. Следовательно, в нашем
  распоряжении оказывается мостик, канонически связывающий свободную симплициальную про-$\mathfrak{C}$-резольвенту с классической
  гомологической теорией групп.

Начнем наши построения с утверждения необходимой терминологии скрещенных модулей и комплексов, следуя \cite{BH,P1,P2}.

\begin{opr}
Про-$\mathfrak{C}$-предскрещенный модуль состоит из про-$\mathfrak{C}$-групп $G_2, G_1$; $G_1$ непрерывно действует слева на $G_2$
$(g_1,g_2) \to {}^{g_1}g_2$,  и непрерывного гомоморфизма про-$\mathfrak{C}$-групп $$\partial: G_2 \to G_1$$
такого, что для всех $g_2 \in G_2, g_1 \in G_1$
$$\text{CM 1). } \partial({}^{g_1}g_2) = g_1 \partial(g_2) g_1^{-1}$$
\end{opr}

\begin{opr}
Про-$\mathfrak{C}$-предскрещенный модуль называется про-$\mathfrak{C}$-скрещенным, если дополнительно выполнено тождество:
$$\text{CM 2). } ^{\partial(g_2)}g_2'= g_2 g_2' g_2^{-1}$$
Тождество CM 2) называется Пайферовым тождеством.

В дальнейшем будем обозначать такой про-$\mathfrak{C}$-скрещенный модуль $(G_2,G_1,\partial)$.
\end{opr}

\begin{opr}
Морфизмом про-$\mathfrak{C}$-скрещенных модулей $(G_2,G_1, \partial) \to (G'_2,G'_1, \partial')$ называется пара непрерывных
гомоморфизмов $$\varphi: G_2 \to G'_2 \text{ и } \psi: G_1 \to G'_1$$ таких, что  $$\varphi(^{g_1}g_2) = ^{\psi(g_1)}\varphi (g_2) \text{
и } \partial' \varphi(g_2) = \psi \partial(g_2).$$
\end{opr}
Возникающая категория про-$\mathfrak{C}$-скрещенных модулей будет обозначаться $ModС^{pro-\mathfrak{C}}$.

Для нас ключевое значение будут играть скрещенные модули, возникающие из копредставлений про-$\mathfrak{C}$-групп. Пусть задано
копредставление $P=(X,R)$ про-$\mathfrak{C}$-группы $G$, т.е. в нашем расположении имеется короткая точная последовательность
$$1 \to N(R) \to F(X) \to G \to 1,$$
где $F$ - свободная про-$\mathfrak{C}$-группа на замкнутом проконечном пространстве $X, R$ - замкнутое подпространство в $F$ и $N=N(R)$ -
замкнутая нормальная подгруппа, порожденная $R$. Группа $F$ действует непрерывно на $N$ сопряжением $g \mapsto ^fg = f g f^{-1}, \, g \in
N, \, f \in F.$

Для дальнейших построений нам потребуется несколько более формальная конструкция, мы рассмотрим копредставление, как тройку $(X, \bar{R},
\omega)$, где $\bar{R}$ - проконечное пространство, $\omega: \bar{R} \to F$ - непрерывное отображение $\bar{R}$ в $F$, где $F$ -
свободная про-$\mathfrak{C}$-группа на проконечном пространстве $X$. Элементы $\bar{R}$ будем обозначать греческими $\rho, \sigma, \tau,
..., $ а элементы $\omega(\rho), \omega(\sigma), \omega(\tau)$,... из $R = \omega(\bar{R})$ латинскими $r,s,t,...$
Пусть $F_\mathfrak{C}(X \cup \bar{R})$ - свободная про-$\mathfrak{C}$-группа на проконечном пространстве $X \cup \bar{R}.$ Пусть $H$ -
это нормальное замыкание $\bar{R}$ в $F_\mathfrak{C}(X \cup \bar{R})$. Построим  $\Theta: H \to F(X),$ который получается ограничением
на $H$ гомоморфизма $\Theta: F_\mathfrak{C}(X \cup \bar{R}) \to F(X),$ заданного на базисных элементах $X \cup \bar{R}$ по формулам:
$\Theta(x) = x, \, x \in X, \, \Theta(\rho) = \omega(\rho) = r, \, \rho \in \bar R, \, \omega(\rho) = r, \, r \in R$.
Ясно, что $\Theta(H) = N(R).$ Если заметить, что $F_\mathfrak{C}(X \cup \bar{R})$ действует на себе слева сопряжением, то $\Theta$
является операторным гомоморфизмом
$\Theta(^fg) = f \Theta(g) f^{-1}, \, f \in F, \, h \in H.$

Очевидно, что для $(F,H,\Theta)$ выполнено $СM 1) \Rightarrow (F,H,\Theta)$ - про-$\mathfrak{C}$-предскрещенный модуль. Элементы $E_P: =
\Ker \Theta$ называются тождествами среди соотношений. Группа $E$ содержит очевидные тождества среди соотношений вида $$(\rho \sigma
\rho^{-1})(^{\rho}\sigma)^{-1}, \text{ где } \rho, \sigma \in \bar{R}$$
Для формализации этого феномена введем Пайферовы коммутаторы $\langle a,b\rangle = ab a^{-1} (^{\Theta(a)}b)^{-1}$
\begin{opr}
Про-$\mathfrak{C}$-скрещенным модулем копредставления называется про-$\mathfrak{C}$-скрещенный модуль $(H/P,F, \partial)$, где $P$ -
Пайферова подгруппа $(H,F, \partial)$, т.е. замкнутая подгруппа в $H$, порожденная Пайферовыми коммутаторами.
\end{opr}
Корректность этого определения немедленно следует из следующего общего предложения.
\begin{predl}\cite{BH}
Пусть $(G_2,G_1,\partial)$ - про-$\mathfrak{C}$-предскрещенный модуль, тогда:

1) Пайферова группа $P$ про-$\mathfrak{C}$-скрещенного модуля - это замкнутая подгруппа $G_2$, порожденная Пайферовыми коммутаторами,
является нормальной и $G_1$-инвариантной.

2) $(G_2/P,G_1,\partial)$ - про-$\mathfrak{C}$-скрещенный модуль и морфизм про-$\mathfrak{C}$-предскрещенных модулей
$$\varphi: (G_2,G_1,\partial) \to (G_2/P, G_1, \partial)$$
обладает универсальным свойством относительно морфизмов из $(G_2,G_1,\partial)$ в про-$\mathfrak{C}$-скрещенные $G_1$-модули.
\end{predl}

Если $G$ - про-$\mathfrak{C}$-группа, а $P=(X,R)$ - её копредставление, то можно дать более точное описание $H$. Пусть $H$ - свободная
$F$-операторная про-$\mathfrak{C}$-группа на пространстве $\bar{R}$ с левым действием $F (h,u)\mapsto ^uh, \, h \in H, \, u \in F,$ т.е.
$H$ - свободная про-р-группа на проконечном пространстве $Y=\bar{R} \times F,$ элементы которого будем записывать в виде $^u\rho, \rho
\in \overline{R}, u \in F, \, ^1\rho$ обозначается, как $\rho, \, ^{-1}(^u\rho):=^{-u}\rho$ или $^u(^{-1}\rho)$. Тогда имеет место
\begin{predl}
\cite[~Corollary 2.2]{Gl1} Пусть $F(X \cup \bar{R})$ - свободная про-$\mathfrak{C}$-группа, порожденная $X \cup \bar{R}$. Тогда $H$
изоморфна нормальному замыканию $\bar{R}$ в $F(X \cup \bar{R})$.
\end{predl}
\textbf{Доказательство:} Отметим, что замкнутость класса $\mathfrak{C}$ относительно расширений существенна при доказательстве
предложения, т.к. используется \cite[теорема 2.1]{Gl1}. $\boxtimes$

Теперь мы можем дать точное описание Пайферовой подгруппы предскрещенного модуля, возникающего из копредставления
про-$\mathfrak{C}$-группы.
\begin{predl}
Пусть $(H,F,\Theta)$ - про-$\mathfrak{C}$-предскрещенный $F$-модуль, соответствующий копредставлению $P=(X,R,\omega)$. Тогда Пайферова
группа $P$ является нормальным замыканием в $H$ базисных Пайферовых элементов $\langle a,b\rangle, \, a,b \in \bar{R} \times F.$
\end{predl}
Доказательство следует из общего утверждения:
\begin{predl}
Пусть $(G_2,G_1,\partial)$ - про-$\mathfrak{C}$-скрещенный $G_1$-модуль, пусть $V$ - множество образующих $G_2$ таких, что $V$ является
$G_1$-инвариантным. Тогда Пайферова подгруппа $P$ совпадает с нормальным замыканием в $G_2$ пространства $Z$ Пайферовых коммутаторов
$$\langle a,b\rangle, \, a,b \in V.$$
\end{predl}
\textbf{Доказательство.} Представляет собой компиляцию рассуждений из дискретного случая  \cite [Предл.3]{BH}.$\boxtimes$

Перечислим основные элементарные свойства про-$\mathfrak{C}$-скрещенных модулей \cite{BH, P1,P2}.

Пусть $(G_2,G_1,\partial)$ - про-$\mathfrak{C}$-скрещенный модуль.
\begin{lemma}["Образы нормальны"]
Пусть $N = \Im \partial$, тогда $N \triangleleft G_1.$
\end{lemma}
\begin{lemma}["Ядра центральны"]
Пусть $A=\Ker \partial.$ Тогда $A$ лежит в центре $G_2$.

\end{lemma}
\begin{lemma}(\textbf{$G_1$-действие превращает $A$ в $\mathbb{\widehat{Z}}_\mathfrak{C}[[G_1]]$-модуль})

Замкнутая нормальная подгруппа $N = \Im \partial \lhd G_1$ действует тривиально на $A$, следовательно $A$ является проконечным левым
$\mathbb{\widehat{Z}}_\mathfrak{C}[[G_1/N]]$-модулем.
\end{lemma}
\begin{lemma}["Абелианизация"]
Абелианизация $G_2^{Ab}=G_2/\overline{[G_2,G_2]}$ имеет структуру левого проконечного
$\mathbb{\widehat{Z}}_\mathfrak{C}[[G_1/N]]$-модуля.
\end{lemma}

Естественно, что и $N^{Ab}$ также имеет структуру проконечного $\mathbb{\widehat{Z}}_\mathfrak{C}[[G_1/N]]$-модуля.

\begin{predl}\cite{BH}
Пусть $(G_2,G_1,\partial)$ - про-$\mathfrak{C}$-скрещенный модуль, тогда индуцированный из $\partial$ морфизм дает точную
последовательность проконечных $\mathbb{\widehat{Z}}_\mathfrak{C}[[G_1/N]]$-модулей
$$A \to G_2^{Ab} \to N^{Ab} \to 0$$
\end{predl}

Вообще говоря, ядро гомоморфизма из $A$ в $G_2^{Ab}$, равное $A \cap \overline{[G_2,G_2]}$ нетривиально. Действительно, следуя
\cite{P1,BH}, рассмотрим скрещенный модуль $\partial: G \to \Aut G$, где $G$ - конечная диэдральная группа, $\Aut G$ действует на $G$
(понятно как), а $\partial$ переводит  элемент $g \in G$ в автоморфизм сопряжения. Ясно, что $\Ker \partial = ZG$ - центр группы, но в
диэдральной группе $ZG \cap [G,G] \ne 1.$

Однако, в случае, когда $\partial : G_2 \to N$ допускает непрерывное расщепление, имеет место
\begin{predl}
Пусть $(G_2,G_1,\partial)$ - про-$\mathfrak{C}$-скрещенный модуль, $1 \to A \to G_2 \to N \to 1$ - точная последовательность $A = \Ker
\partial, \, N = \Im \partial$ и эпиморфизм имеет непрерывное расщепление, т.е. существует $s: N \to G_2: \partial \circ s = \id_N,$ то
имеет место точная последовательность
$$0 \to A \to G_2^{Ab} \to N^{Ab} \to 0$$
\end{predl}
\textbf{Доказательство.}
Пусть задано непрерывное расщепление $s: N \to G_2,$ тогда $G_2=A \la s(N).$ Т.к. $A$ - центральна, то все коммутаторы лежат в $s(N)$, а,
следовательно,
$A \cap \overline{[G_2,G_2]} \subseteq A \cap s(N) = 1. \boxtimes$

Важным классом про-$\mathfrak{C}$-скрещенных модулей, к которым относятся, в частности, про-$\mathfrak{C}$-скрещенные модули, возникающие
из копредставлений про-$\mathfrak{C}$-групп, являются свободные про-$\mathfrak{C}$-скрещенные модули. Их особенностью является
существование непрерывного расщепления, а следовательно и точность последовательности $0 \to A \to G_2^{Ab} \to N^{Ab} \to 0,$ что станет
ключевым моментом в реализации нашей задачи по установлению связи между скрещенными модулями копредставлений групп и последовательностью
Кроуэла-Линдона, которая определяет начало канонической $\mathbb{\widehat{Z}}_\mathfrak{C}[[G]]$-резольвенты копредставления
про-$\mathfrak{C}$-группы $G$.

Пусть $(G_2,G_1,\partial)$ - про-$\mathfrak{C}$-предскрещенный модуль, $\bar{R}$ - проконечное пространство и $f: \bar{R} \to G_2$ -
непрерывное отображение. Тогда $(G_2,G_1,\partial)$ - свободный про-$\mathfrak{C}$-предскрещенный модуль с \textbf{базисом $f$}, если для
любого про-$\mathfrak{C}$-предскрещенного $G_1$-модуля $(G'_2,G_1,\partial')$ и функции $f': \overline{R} \to G'_2$ такой, что $\partial'
f' = \partial f$, то существует и единственный морфизм про-$\mathfrak{C}$-предскрещенных модулей $\varphi:(G_2,G_1,\partial) \to
(G'_2,G_1,\partial'): \varphi f = f'.$ Для про-$\mathfrak{C}$-скрещенного модуля, возникающего из копредставления, естественно
рассматривать также и непрерывное отображение $\omega= \partial f: \overline{R} \to G_1,$ называя такой про-$\mathfrak{C}$-скрещенный
модуль $(G_1,G_2,\partial)$ с функцией $f$ свободным про-$\mathfrak{C}$-предскрещенным модулем на $\omega$.

Если $(G_2,G_1,\partial)$ - скрещенный про-$\mathfrak{C}$-модуль, который обладает описанным выше универсальным свойством в категории
про-$\mathfrak{C}$-скрещенных модулей, то мы назовем $(G_2,G_1,\partial)$ свободным про-$\mathfrak{C}$-скрещенным модулем с базисом $f$
(или на $\omega$).
\begin{predl} \cite{BH}
Пусть $G$ - про-$\mathfrak{C}$-группа, $\overline{R}$ - проконечное пространство, и $\omega: \overline{R} \to G$ - непрерывное
отображение. Тогда существуют и единственны

(i) свободный про-$\mathfrak{C}$-предскрещенный модуль на $\omega$

(ii) свободный про-$\mathfrak{C}$-скрещенный модуль на $\omega$
\end{predl}

Из конструкции свободного про-$\mathfrak{C}$-предскрещенного модуля на $\omega$ ясно, что базисное отображение $f: \bar{R} \to H$
инъективно. Оно будет инъективно и для про-$\mathfrak{C}$-скрещенного модуля. Это очевидно следует  из следующего более сильного
утверждения
\begin{predl}\cite{BH}
Пусть $(G_2,G_1,\partial)$ - свободный про-$\mathfrak{C}$-скрещенный модуль с базисом $f:\bar{R}\to G_2, \bar{R}$ - проконечное
пространство, $G = Coker \partial,$ тогда $G_2^{Ab}$ является свободным проконечным $\widehat{\mathbb{Z}}_\mathfrak{C}[[G]]$-модулем на
композиции $\bar{f}: \bar{R} \rightarrow G_2 \to G_2^{Ab}$ и, следовательно, $\bar{f}$ (тогда и $f$) - инъективные отображения.
\end{predl}

\begin{theorem}\cite{BH}
Пусть $(G_2,G_1,\partial)$ - свободный про-$\mathfrak{C}$-скрещенный модуль, построенный по копредставлению $P=(X,\bar{R},\omega)$
некоторой про-$\mathfrak{C}$-группы $G$. Пусть $\pi_P = \Ker \partial$ и $N = \Im \partial$ - свободная про-$\mathfrak{C}$-группа (в
про-р-случае $N$ всегда свободна, как подгруппа свободной про-р-группы). Тогда индуцированный гомоморфизм $\xi: \pi_P \to G_2^{Ab}$
инъективен и возникает короткая точная последовательность $G$-модулей $(G \cong G_1/N)$
$$0 \to \pi_P \to G_2^{Ab} \xrightarrow{d} N^{Ab} \to 0,$$ в которой $G_2^{Ab}$ является свободным
$\mathbb{\widehat{Z}}_\mathfrak{C}[[G]]$-модулем на элементах $\bar{f}(\rho) = \rho\overline{[G_2,G_2]}, \, \rho \in \bar{R},$ а $d$
задано формулой $d(\rho\overline{[G_2,G_2]}) = \omega(\rho)[N,N]$ \label{20}
\end{theorem}

 $N^{Ab}$ - в литературе фигурирует, как модуль соотношений, оказывается, что $N^{Ab}$ естественным образом вкладывается в
 $\mathbb{\widehat{Z}}_\mathfrak{C}[[G]]$-модуль $\mathbb{\widehat{Z}}_\mathfrak{C}[[G]]^{(X)}$.
К деталям этого вложения мы вернемся несколько позже, а пока покажем, как из свободной про-$\mathfrak{C}$-симплициальной резольвенты,
возникающей из некоторого копредставления про-$\mathfrak{C}$-группы $G$, получается про-$\mathfrak{C}$-скрещенный модуль, а также, что
этот модуль совпадает с введенным ранее свободным про-$\mathfrak{C}$-скрещенным модулем копредставления про-$\mathfrak{C}$-группы.

Пусть задано копредставление $P=(X,\bar{R},\omega)$ про-$\mathfrak{C}$-группы $G$, это служило исходным материалом для построения
свободной симплициальной про-$\mathfrak{C}$-резольвенты в первой части работы. Имеет место короткая точная последовательность
$1 \to N(\bar{R}) \to F_\mathfrak{C}(X) \to G \to 1,$
элементы $\omega(\rho)$ порождают $N(\bar{R})$, как топологический нормальный делитель, возникают нулевой и первый шаги построения
свободной симплициальной про-$\mathfrak{C}$-резольвенты
$
\xymatrix{F_\mathfrak{C}(X_0 \cup \bar{R}) \ar@<1ex>[r]^{d_1^1} \ar@<-1ex>[r]_{d_0^1} & F_\mathfrak{C}(X_0) \ar@{->>}[r]^{\varepsilon} &
G
},$
где $d_0^1(x) = x, \, x \in X_0, \, d_0^1(\rho) = 1, \, \rho \in \bar{R},$
$d_1^1(x)=x, \, x \in X_0, \, d_1^1(\rho) = \omega(\rho).$

Метод "pas a pas" даёт на первом шаге свободную симплициальную про-$\mathfrak{C}$-группу $F^{(1)}$, у которой $F_n^{(1)}: F_n^{(1)} \cong
F_\mathfrak{C}(\bigcup_{n \twoheadrightarrow [1]} s_i(X_1) \cup X_0)$ при  $n \ge 2$ вырождены.

Пусть $\partial = d_1^1: NF_1^{(1)} \to NF_0^{(1)} = F_0^{(1)}$ ($NF^{(1)}$- комплекс Мура симплициальной про-$\mathfrak{C}$-группы
$F^{(1)}$).  $F_0^{(1)}$ действует слева на $NF_1^{(1)}$ сопряжением, т.е. $^gx = s_0(g) x s_0(g)^{-1}.$
Тогда условие $CM\,$1) очевидно выполнено, т.к. $\partial(^gx) = \partial (s_0(g)xs_0(g)^{-1}) = g \partial(x) g^{-1}.$

Следовательно $(NF_1^{(1)}, N F_0^{(1)} =F_0^{(1)}, \partial)$ является про-$\mathfrak{C}$-предскрещенным модулем.

При построении про-$\mathfrak{C}$-скрещенного модуля копредставления, в общем случае, условие $CM\,$2)
$x y x^{-1} = ^{\partial x}y = ^{d_1^1(x)}y$
конечно же выполнено не всегда. Сформируем Пайферовы коммутаторы
\begin{equation}
\langle x,y\rangle = xyx^{-1}(^{d_1^1(x)}y)^{-1}=xyx^{-1}(s_0 d_1^1(x) y s_0 d_1^1(x)^{-1})^{-1}
\end{equation}

Если произвести замену $x\mapsto x^{-1}, \, y\mapsto x^{-1}y^{-1}x$, то получим, что $\langle x,y\rangle \approx \langle
x,y\rangle^{-1}=[x^{-1}s_0d_1^1(x),y]$. Следовательно, Пайферова подгруппа порождается элементами вида $[x^{-1} s_0 d_1^1(x), y].$

Теперь заметим, что любой элемент $x \in \Ker d_1^1$ может быть записан в виде $s_1 d_1(y) \cdot y^{-1},$ где $y \in NG_1=\Ker d_0^1,$
действительно, т.к.
$G_1 \cong NG_1 \la s_0 G_0,$ то $x=y \cdot s_0(y_0),$ где $y_0 \in G_0, \, y \in NG_1.$
Т.к. $x \in \Ker d_1^1$, имеем
$1 = d_1^1(x) = d_1^1(y) \cdot d_1^1 s_0(y_0) \Rightarrow y_0 = d_1^1(y^{-1}) \Rightarrow x=y \cdot s_0 d_1^1(y^{-1})$
Следовательно, Пайферова подгруппа совпадает с $[\Ker d_1^1; \Ker d_0^1],$ откуда получаем, что про-$\mathfrak{C}$-скрещенный модуль
изоморфен $(\frac{\Ker d_0^1 }{[\Ker d_0^1, \Ker d_1^1]}, F_0^{(1)}, d_1^1).$

Заметим, что
$\hr{\frac{NF_1^{(1)}}{d_2^2 NF_2^{(1)}},NF_0^{(1)} = F_0, \partial=d^1_1}$
является скрещенным модулем. Действительно, $d_2^2 N F_2^{(1)}$ - замкнутая нормальная подгруппа в $ \Ker d_1^1,$ поэтому $\partial =
d_1^1$ определен корректно. Действие $NF_0^{(1)}$ индуцировано действием $F_0$ путем сопряжения на элемент $s_0(g)$ тоже корректно, т.к.
$d_2^2 NF_2^{(1)}$ является $F_0^{(1)}$-инвариантной подгруппой.
Выполнение $CM\,$1) следует из выполнения $CM\,$1) для про-$\mathfrak{C}$-предскрещенного модуля $NF_1^{(1)}\xrightarrow{\partial =
d_1^1} NF_0^{(1)}.$ Остается проверить выполнение  $CM\,$2), т.е. надо проверить  в $\frac{NF_1^{(1)}}{d_2^2 NF_2^{(1)}}$ равенство
$\langle x,y\rangle = 1.$
Действительно, $\langle x,y\rangle \cong (s_0 d_1^1(x^{-1}) y s_0 d_1^1(x))(x^{-1} y x)^{-1}\cong [s_0 d_1^1(x), y^{-1}][y^{-1},x] \cong
[s_0d_1^1(x),y][y,x].$
Рассмотрим элемент
$\{x,y\}:=[s_0x,s_1y][s_1y,s_1x],$
который называется Пайферовским поднятием. Элементарная проверка дает
$d_0^2\{x,y\} = d_1^2\{x,y\} = 1,\, d_2^2\{x,y\}=\langle x,y\rangle.$
Следовательно $\langle x,y\rangle = 1$ в группе $\frac{NF_1^{(1)}}{\partial_2 NF_2^{(1)}}.$

\begin{predl}
Пусть $F^{(1)}$ первый шаг построения симплициальной про-$\mathfrak{C}$-резольвенты методом "pas a pas", тогда скрещенный
$F_0^{(1)}$-модуль $\frac{NF_2^{(1)}}{d_2^2 NF_2^{(1)}} \xrightarrow{d_1^1} F_0^{(1)}$ совпадает с про-$\mathfrak{C}$-скрещенным
$F$-модулем, полученным из копредставления $(H/P,F,\partial)$.

\end{predl}

\textbf{Доказательство.}
Для доказательства нам потребуется лемма Брауна-Лодэя, полезная в приложениях симплициальной теории групп.
\begin{lemma}\cite{P2}
Пусть $G$ - симплициальная про-$\mathfrak{C}$-группа такая, что $G_2$ вырождена, т.е. $G_2$ = $D_2= \{$подгруппа, порождённая $s_0(G_1),
s_1(G_1)\}$, тогда замкнутая подгруппа, порожденная коммутаторами $[x,y],$ где $x \in \Ker d_0^1, y \in \Ker d_1^1$ совпадает с $d_2^2 N
G_2: \overline{[\Ker d_0^1, \Ker d_1^1]} = d_2^2 N G_2.$ \label{21}
\end{lemma}

\begin{zamech}
Доказательство леммы \ref{21} можно найти в \cite{Br2},\cite[~prop.23]{P1} для про-$\mathfrak{C}$-случая, оно опирается на разложение
Кондуше:
\begin{utv}\cite[~стр. 158]{Con}
Если $G$ - симплициальная про-$\mathfrak{C}$-группа, то для любого $n \ge 0$ справедливо
$$G_n \cong (...(N G_n \times s_{n-1} N G_{n-1}) \times ... \times s_{n-2}...s_1 N G_1) \times (...(s_0 N G_{n-1} \times s_1 s_0 N
G_{n-2}) \times ... \times s_{n-1} s_{n-2}...s_0 N G_0).$$

Расстановка скобок и порядок составляющих этого полупрямого произведения порождается последовательностью:
$$G_1 \cong N G_1 \leftthreetimes s_0 N G_0,$$
$$G_2 \cong (N G_2 \leftthreetimes s_1 N G_1)\leftthreetimes (s_0 N G_1\leftthreetimes s_1 s_0 N G_0),$$
$$G_3 \cong ((N G_3\leftthreetimes s_2 N G_2)\leftthreetimes (s_1 N G_2\leftthreetimes s_2 s_1 N G_1))\leftthreetimes ((s_0 N
G_2\leftthreetimes s_2 s_0 N G_1)\leftthreetimes (s_1 s_0 N G_1\leftthreetimes s_2 s_1 s_0 N G_0)),$$
и т.д.
\end{utv}

Разложение Кондуше $G_2 \cong (N G_2\leftthreetimes s_0 N G_1)\leftthreetimes(s_1 N G_1\leftthreetimes s_1 s_0 N G_0)$ вместе с индукцией
по длине выражения, представляющего элемент из $N G_2$, приводит к доказательству Леммы \ref{21}.
\end{zamech}

Уже доказано, что $(H/P,F,\partial) \cong \hr{\frac{N F_1^{(1)}}{\overline{[\Ker d_0^1, \Ker d_1^1]}}, F_0^{(1)}, d^1_1}$, т.к.
$F_2^{(1)} = D_2$($F_1^{(1)})$ из леммы Брауна-Лодэя следует, что $\hr{\frac{N F_1^{(1)}}{\overline{[\Ker d_0^1, \Ker d_1^1]}},
F_0^{(1)}, d_1^1} \cong \hr{\frac{N F_1^{(1)}}{d_2^2 N F_2^{(1)}},
F_0^{(1)}, d_1^1}. \boxtimes$

Имеем точную последовательность
$\pi  = \pi_1(F^{(1)}) \to \frac{NF_1^{(1)}}{d_2^2 NF_2^{(1)}} \xrightarrow{d_1^1} \Im d_1^1$.
Из теоремы \ref{20} получаем
\begin{sled}
$\pi_1(F^{(1)}) \cong \pi_P(G)$, где $\pi_P(G)$-модуль тождеств среди соотношений.
\end{sled}

\begin{opr}
Про-$\mathfrak{C}$-скрещенный комплекс $C$ состоит из про-$\mathfrak{C}$-групп и непрерывных морфизмов
$$C: ... \to C_n \xrightarrow{\delta_n} C_{n-1} \xrightarrow{\delta_{n-1}} ... C_2 \xrightarrow{\delta_2} C_1 \xrightarrow{\delta_1}
C_0,$$
удовлетворяющих следующим условиям:

СС 1) $\delta_1: C_1 \to C_0$ - про-$\mathfrak{C}$-скрещенный модуль;

СС 2) $C_n(n >1)$ является проконечным левым $\widehat{\mathbb{Z}}_\mathfrak{C}[[C_0/\delta_1 C]]$-модулем и каждый $\delta_n(n >1)$ -
непрерывный гомоморфизм левых проконечных $\widehat{\mathbb{Z}}_\mathfrak{C}[[C_0/\delta_1 C]]$-модулей (для $n=~2$ $ \delta_2$
коммутирует с действиями $C_0$ и $\delta_2(C_2) \subset C_1,$ является проконечным $\widehat{\mathbb{Z}}_\mathfrak{C}[[C_0/\delta_1
C]]$-модулем);

СС 3) $\delta \delta = 0.$
\end{opr}
Морфизмы скрещенных комплексов определяются естественным образом.
\begin{opr}
Про-$\mathfrak{C}$-скрещенная резольвента про-$\mathfrak{C}$-группы $G$ - это  про-$\mathfrak{C}$-скрещенный комплекс $C$, где $C_n$ -
проективные $\widehat{\mathbb{Z}}_\mathfrak{C}[[G]]$-модули при $n > 1,$ и выполнено условие точности $ \Im \delta_n = \Ker
\delta_{n-1}$, при этом существует непрерывный изоморфизм $C_0/\delta_1 C_1 \cong G.$
\end{opr}
Про-$\mathfrak{C}$-скрещенные резольвенты могут быть построены из про-$\mathfrak{C}$-копредставлений $P(X,\overline{R},\omega)$.

В качестве скрещенного модуля $C_1 \xrightarrow{\partial} C_0$ берем скрещенный модуль копредставления, т.е. $C_0 = F_\mathfrak{C}(X),$ а
$C_1 = H/P,$ где $H$ - это нормальное замыкание $\overline{R}$ в $F_\mathfrak{C}(X \cup \overline{R})$, а $P$ - Пайферова подгруппа в $H$
и пусть $\pi_P:=\Ker (\partial: C_1 \to C_0)$
- модуль тождеств среди соотношений. Мы ранее доказали, что он является левым проконечным
$\widehat{\mathbb{Z}}_\mathfrak{C}[[G]]$-модулем. Т.к. категория проконечных модулей имеет достаточно проективных объектов \cite{ZR}, мы
можем построить $\widehat{\mathbb{Z}}_\mathfrak{C}[[G]]$-проективную резольвенту $\widehat{\mathbb{Z}}_\mathfrak{C}[[G]]$-модуля $\pi_P.$
Для получения скрещенной резольвенты $C$ про-$\mathfrak{C}$-группы $G$ положим $C_n = P_{n-2}, \, n > 2, \, \delta_n = d_{n-2}, \, n >2,$
и  $d_2$ - композит $P_0 \to \pi_P \xrightarrow{i} C_1.$

Касаемо примеров про-$\mathfrak{C}$-скрещенных резольвент и их свойств, см. \cite{P1,P2}. В частности имеется стандартная неоднородная
про-$\mathfrak{C}$-резольвента, являющаяся прямым аналогом неоднородной бар-резольвенты.

Из теоремы \ref{20} следует, что $C_1^{Ab} \cong \mathbb{\widehat{Z}}_C[[G]]^{(R)}$ - свободный проконечный
$\widehat{\mathbb{Z}}_\mathfrak{C}[[G]]$-модуль на $R$, a $\partial_{ab}$ индуцирует гомоморфизм $\partial_*:
\mathbb{\widehat{Z}}_\mathfrak{C}[[G]]^{(R)} \to N(R)^{Ab},$ заданный формулой $\partial_*(e_{\rho}) = r[\overline{N(R),N(R)}],$ где
$e_{\rho}$ - образующая $\widehat{\mathbb{Z}}_\mathfrak{C}[[G]]^{(R)},$ соответствующая элементу $\rho \in \bar{R}$.
С другой стороны, точной последовательности копредставления $$1 \to N(R) \xrightarrow{i} F_\mathfrak{C}(X) \xrightarrow{\varphi} G \to
1$$ соответствует \cite[теор.~6.3]{HS} короткая точная последовательность

$$0 \to N(R)^{Ab} \xrightarrow{\bar{i}} \mathbb{\widehat{Z}}_\mathfrak{C}[[G]]\widehat{\otimes}_F I_\mathfrak{C}(F)
\xrightarrow{\bar{\varphi}} I_\mathfrak{C}(G) \to 0,$$
где $F = F_\mathfrak{C}(X).$ С другой стороны, $I_\mathfrak{C}(F)$ является свободным модулем, порожденным элементами $(1-x)$, где $x \in
X.$ Получаем
$\widehat{\mathbb{Z}}_\mathfrak{C}[[G]] \widehat{\otimes}_F I_\mathfrak{C}(F) \cong \widehat{\mathbb{Z}}_\mathfrak{C}[[G]]
\widehat{\otimes}_F \mathbb{\widehat{Z}_\mathfrak{C}}[[G]]^{(X)} \cong \mathbb{\widehat{Z}_\mathfrak{C}}[[G]]^{(X)}$.
Откуда $$\widehat{\mathbb{Z}}_\mathfrak{C}[[G]]^{(R)} \cong C_1^{Ab} \xrightarrow{\partial_*} N(R)^{Ab} \xrightarrow{i}
\mathbb{\widehat{Z}}_\mathfrak{C}[[G]]\widehat{\otimes}_F I_\mathfrak{C}(F) \to \mathbb{\widehat{Z}}_\mathfrak{C}[[G]]^{(X)}$$

Если заменить $I_\mathfrak{C}(G)$ на $\Ker(\mathbb{\widehat{Z}}_\mathfrak{C}[[G]] \to \widehat{\mathbb{Z}}_\mathfrak{C}),$ то получим
свободную проконечную $\mathbb{\widehat{Z}}_\mathfrak{C}[[G]]$-резольвенту тривиального модуля $\widehat{\mathbb{Z}}_\mathfrak{C}$
$$\mathbb{\widehat{Z}}_\mathfrak{C}[[G]]^{(R)} \xrightarrow{d} \mathbb{\widehat{Z}}_\mathfrak{C}[[G]]^{(X)} \to
\mathbb{\widehat{Z}}_\mathfrak{C}[[G]] \to \mathbb{\widehat{Z}}_\mathfrak{C}, d = i \circ \partial_*,$$
построенную из копредставления, которая является про-$\mathfrak{C}$-аналогом последовательности Кроуэла-Линдона комплекса цепей
универсального накрытия $\widetilde{K(G,1)},$ где $K(G,1)$ строится по копредставлению стандартным образом, путем приклеивания клеток,
соответствующих соотношениям в $P$ \cite{BH}.

В дискретной ситуации гомоморфизм $\mathbb{Z}[G]$-модулей $d$ допускает удобное описание на языке свободного дифференциального
исчисления, а именно, если $e^2_{\rho}$ - образующий в $Z[G]^{(R)}$, соответствующий соотношению $\rho$, имеющему в $F$ запись
$\omega(\rho) = y_1...y_m,$ где $y_i \in X \cup X^{-1}, i = 1, ...,n,$ то
$d(e_{\rho}^2) = \sum_{x \in X} \varphi (\frac{\partial r}{\partial x}) e_x ,$
где $e_x$ - образующие $\mathbb{Z}[G]$-модуля $\mathbb{Z}[G]^{(X)},$ соответствующие элементам $1-x$ в $I(X),$ $\frac{\partial
r}{\partial x}$ - производные Фокса элемента $r$, $\varphi: F \twoheadrightarrow G$ - каноническая сюръекция. Что касается производных
Фокса $\frac{\partial}{\partial x}: F \to \mathbb{Z}[F] $ для каждого $x \in X$, то они определяются
$$ (i)\quad \frac{\partial y}{\partial x}= \begin{cases}
1, & x = y\\
0, & x \ne y
\end{cases}, y \in X; \quad(ii) \text{ для } \omega_1 \cdot \omega_2 \in F \quad \frac{\partial}{\partial x}(\omega_1\omega_2) =
\frac{\partial}{\partial x}(\omega_1) + \omega_1\frac{\partial}{\partial x}(\omega_2). $$

Вообще говоря, не существует аналогичного описания $d$ в категории $\widehat{\mathbb{Z}}_\mathfrak{C}[[G]]$-модулей. Однако, в своей
работе \cite{I}, Ихара построил аналогичную дискретному случаю теорию свободного дифференциального исчисления для свободных почти
про-р-групп.

Пусть $F$ - абстрактная свободная группа ранга $r$, порожденная элементами $x_1,...,x_r (r \ge 1), F_1 \triangleleft F$ нормальная
подгруппа конечного индекса, тогда $F_1 $ - свободная группа ранга $r_1$ и $r_1 - 1= (r-1)(F:F_1).$ Пусть $\widehat{N}$ - подгруппа
конечного индекса в $F$, тогда обозначим $ N = \widehat{N} \cap F_1$ нормальную подгруппу конечного индекса в $F_1$. Рассмотрим такие
$\widehat{N} \triangleleft F$, что $F_1/ N$ является конечной $p$-группой. Сформируем проективный предел
$\mathfrak{F} = \varprojlim_{\begin{matrix} \widehat{N} \lhd_0 F; [F_1:N]=p^{\alpha}, \\ N = \widehat{N}\cap  F_1 \end{matrix}} F/N.$
Будем называть проконечную группу $\mathfrak{F}$ свободной почти про-р-группой ранга $r$. Она содержит открытую подгруппу $\mathfrak{F_1}
= \varprojlim F_1/N,$ которая является свободной про-р-группой "корректного ранга" $r_1: r_1-1 = (\mathfrak{F}:\mathfrak{F})(r-1).$

В обратную сторону, пусть $\mathfrak{F}$ - проконечная группа, порожденная $r$ элементами, содержит открытую нормальную подгруппу
$\mathfrak{F}_1$, которая является свободной про-р-группой "корректного ранга", тогда мы утверждаем, что $\mathfrak{F}$ является
свободной почти про-р-группой ранга $r$ относительно любого множества из $r$ образующих $\bar{x}_1,...,\bar{x}_r$ в $\mathfrak{F}$.
Действительно, рассмотрим гомоморфизм $\varphi$ из свободной дискретной группы $F=\langle x_1,...,x_r\rangle$ ранга $r$ в $\mathfrak{F}$,
определенный по правилу $x_j \mapsto \bar{x}_j (1 \le j \le r).$ Этот гомоморфизм будет инъективным. На самом деле, пусть $F_1 =
\varphi^{-1}(\mathfrak{F}_1)$.
Тогда ввиду того, что $\overline{\varphi(F)} = \mathfrak{F}$ получаем, что $[F:F_1] = [\mathfrak{F}:\mathfrak{F}_1]$. Следовательно, т.к.
$F$ - свободна, по теореме Шрайера получаем, что $\rk F_1 = \rk \mathfrak{F}_1.$ Рассмотрим свободную почти про-р-группу $\widehat{F}:$
$\widehat{F} \cong \varprojlim_{{\small\begin{matrix} \widehat{N} \lhd_0 F; [F_1:N]=p^{\alpha} \\ N = \widehat{N} \cap F_1 \end{matrix}}}
F/N$
и $\bar{\varphi}$-гомоморфизм, индуцированный $\varphi$.

Мы утверждаем, что $\overline{\varphi} : \widehat{F} \xrightarrow{\bar{\varphi}} \mathfrak{F}$ является изоморфизмом, при этом
$\overline{\varphi} : \widehat{F}_1 \xrightarrow{\bar{\varphi}} \mathfrak{F_1}.$ Ясно, что $\bar{\varphi}$ является сюръекцией. Можно
рассмотреть конфинальную систему нормальных делителей $\mathcal{N}$,в которой все $N$ содержатся в $F_1$, тогда для проверки
инъективности $\bar{\varphi}$ достаточно проверить, что $\bar{\varphi}|_{\widehat{F}_1}$ инъективно, но это следует из известного факта о
том, что свободная про-р-группа конечного ранга не может быть изоморфна своему эпиморфному образу. Ясно, что если свободная почти
про-р-группа является про-р-группой, то она является свободной про-р-группой; также ясно, что любая открытая подгруппа свободной почти
про-р-группы конечного ранга является снова свободной почти про-р-группой корректного ранга. В дальнейшем мы будем отождествлять $x_j \in
F$ с его образом $\bar{x}_j \in \mathfrak{F} (1 \le j \le r).$

\begin{theorem} \cite{I}
Пусть $F$ - свободная почти про-р-группа ранга $r$, порожденная $x_1,...,x_r,$ $B = \mathbb{Z}_p[[\mathfrak{F}]]$ - пополненная групповая
алгебра $\mathfrak{F}$ над $\mathbb{Z}_p$ и $s:B \to \mathbb{Z}_p$ - гомоморфизм аугментации. Тогда каждый элемент $\Theta \in B$, может
быть представлен единственным образом в форме:
$$\Theta = s(\Theta) \cdot 1   + \sum_{j=1}^r \Theta_j (x_j -1) , \quad (\Theta_1, ..., \Theta_r \in B),$$ где $1 = 1_\mathfrak{F}$ -
единица $\mathfrak{F}$. Более того для каждого $j$ соответствия $\Theta \to \Theta_j$ определяют непрерывные $\mathbb{Z}_p$-линейные
отображения $B$ в себя, которые будут обозначаться $\frac{\partial \Theta}{\partial x_j} = \Theta_j (1 \le j \le r ).$
\end{theorem}

Из определения немедленно следуют все стандартные правила дифференцирования \cite{I}.

\begin{prim}[Дискретные группы сквозь р-адические очки]

"Теология"(в терминологии Дж.Адамса \cite{A}) пополнений и локализаций позволяет взглянуть по-новому даже на хорошо изученные объекты.
Адамс замечает, что "..Локализацию можно использовать для построения конечных H-пространств, отличных от классических. Например, можно
построить конечный комплекс X, который является H-пространством и на 2-адический взгляд похож на Sp(2), а на 3-адический на $S^3 \times
S^7$. Это пространство X - и не Sp(2) (посмотрите на него сквозь 3-адические очки) и на $S^3 \times S^7$ (рассмотрите его с 2-адической
позиции)... Он чем-то напоминает монстров из средневековых сказок: голова льва на теле лошади. Ещё более жуткий способ заключается в том,
чтобы взять две части животного и составить их, повернув одну из частей: возьмите льва, отрубите ему голову, а затем приделайте её,
обратив пастью назад. Именно так строится пример Хилтона-Ройтберга; роль льва играет пространство Sp(2).."

Примерно та же ситуация возникает при рассмотрении про-р-пополнений дискретных групп, их свойства могут катастрофически отличаться от
своих дискретных "родителей" уже в случае групп с одним соотношением.
Действительно, пусть G задана копредставлением
$P = (x_1,x_2:x_1^l[x_2,x_1^l])$, тогда, т.к. соотношение не является степенью \cite{Gl2,DV}, по этому копредставлению можно построить
K(G,1), который будет 2-мерным асферическим пространством. Однако сквозь p-адические очки мы увидим совершенно другое "существо".
Остановимся на этой группе немного подробнее. Если p не равно l, то соотношение не принадлежит $F^p[F,F]$, следовательно $\widehat{G^p}
\cong F(x_2)$, т.е. сквозь эти очки мы увидим обычную окружность. Оденем теперь l-адические очки, т.е. l=p.

В качестве примера вычислений производных Фокса, вычислим
$\quad d(e_r) = \varphi \hr{\frac{\partial r}{\partial x_1}} e_1 + \varphi \hr{\frac{\partial r}{\partial x_2}}e_2$
Элемент $\al \in \mathbb{Z}_p[[G]]$ будет лежать в $\Ker d \, \Leftrightarrow$
$\al \cdot \varphi\hr{\frac{\partial r}{\partial x_1}} = 0 \quad \text{и} \quad \al \cdot \varphi\hr{\frac{\partial r}{\partial x_2}} =
0$
$$\frac{\partial r}{\partial x_1} = \sum_{k=0}^{p-1} x_1^k+ x_1^p x_2^{-1} x_1^{-p}(x_2 - 1) \cdot \sum_{k=0}^{p-1} x_1^k,
\quad \frac{\partial r}{\partial x_2} = x_1^p x_2^{-1} (x_1^{-p} - 1).$$
Пусть теперь $g_1 = \varphi(x_1), g_2 = \varphi(x_2),$ тогда в $G$ имеем
\begin{equation}
g_1^p = [g_1^p,g_2].\label{9}
\end{equation}
Если $\gamma_k(G)$ - нижний центральный ряд в $G$. Из \ref{9} видно, что $g_1^p \in \gamma_{k+1}(G)$, если $g_1^p \in \gamma_k(G),$ то
$g_1^p \in \gamma_{k+1}(G)$, но $\bigcap_{k=0}^{\infty}\gamma_k(G) = \{1\} \Rightarrow g_1^p = 1 \Rightarrow g_1$ либо имеет порядок $p$,
либо равен 1 в группе. Однако $x_1 \notin N(r),$ т.к. $N(r) \subseteq F^p[F,F] \Rightarrow g_1^p =1,$ тогда
$\varphi \hr{\frac{\partial r}{\partial x_1}} = (2 - g_2^{-1})(1 + g_1+ ...+g_1^{p-1}), \varphi\hr{\frac{\partial r}{\partial x_2}}=0.$
Т.к. $(1-g_1)$ аннулирует $\varphi\hr{\frac{\partial r}{\partial x_1}},$ то
$K(P) \supseteq (1-g_1) \mathbb{Z}_p[[G]] e_r$.

Что касается установления равенства $K(P) = (1-g_1)\mathbb{Z}_p[[G]]e_r,$ то это проще всего установить следующим образом (набросок
построений). Т.к. $g_1^p = 1$ в $G$, то $x_1^p \in R,$ тогда $x_1^p[x_2,x_1^p] = x_1^p x_2^{-1} x_1^{-p} x_2 x_1^p = x_1^p \mod
R^p[R,F]$. Следовательно \cite{Se1}, $x_1^p$ может быть выбран в качестве определяющего соотношения $G$, тогда $G$ может быть
представлена в виде свободного произведения циклической группы $C$ порядка p и свободной про-р-группы ранга 1. Дальнейшие аргументы
понятны, т.к. $I(G)\cong Z_p[[G]]\oplus Z_p[[G]]\widehat{\otimes}_{Z_pC} IC$, \cite[лемма 3.8]{Qui5} см. также \cite{Ku}.
В частности, получим периодическую резольвенту
$$\to \mathbb{Z}_p[[G]] \xrightarrow{1 - g_1} \mathbb{Z}_p[[G]] \xrightarrow{N} \mathbb{Z}_p[[G]] \xrightarrow{1 - g_1}
\mathbb{Z}_p[[G]]^{(1)} \xrightarrow{\partial} \mathbb{Z}_p[[G]]^{(2)} \to \mathbb{Z}_p[[G]] \to \mathbb{Z}_p,\,
N:=1+g_1+...+g_1^{p-1}$$

Эту резольвенту следует воспринимать, как тензорно умноженную на $Z_p[[G]]$ резольвенту циклической подгруппы $C_p < G,$ порожденной
элементом $\langle g_1\rangle,$ которую приклеили к ядру точной последовательности Кроуэла-Линдона $\mathbb{Z}_p[[G]]^{(1)}
\xrightarrow{\partial} \mathbb{Z}_p[[G]]^{(2)} \to \mathbb{Z}_p[[G]] \to \mathbb{Z}_p$. Таким образом, в p=l-адических очках мы видим
следы (на уровне цепного комплекса) линзового пространства.

Этот пример скорее правило, чем исключение. Мельников \cite{Mel2} показал, что для любого элемента, являющегося $p$-ой степенью в
свободной про-р-группе, можно построить элемент $\overline{r}$, который уже не является $p$-ой степенью, однако $r$ и $\overline{r}$
порождают один и тот же нормальный делитель в группе $F$.

До сих пор остается открытым вопрос Серра \cite{Se3}, слегка подправленный Гильденхьюзом, а именно: Верно ли, что если $G$ - про-р-группа
с одним соотношением с $cd(G) > 2$, то можно найти такое копредставление $F/(r)$ с $r = u^p$. Некоторые частичные результаты по этому
вопросу можно найти в \cite{L}.
Заметим \cite{DV}, что K(G,1) для любой дискретной группы с одним соотношением получается либо асферическим, либо путём приклеивания
линзового пространства к $K^1$ по соотношению. Тогда, если подключить фантазию, вопрос о структуре про-р-групп с одним соотношением суть
утверждение о том, что сквозь p-адические очки в дискретных группах с одним соотношением ничего иного, чем линзовые пространства,
асферические комплексы, либо "цветы"(свободные про-р-группы), увидеть не получится.

\end{prim}

Обращают на себя внимание несколько интересных фактов: во-первых, утверждение \cite{GK} о том, что дефицит конечной группы не зависит от
копредставления в категории проконечных групп; во-вторых, согласно \cite[Предложение 7.7]{KH}, \cite{Mel1} свойство про-р-группы быть
асферичной не зависит от копредставления в категории про-р-групп. Таким образом, проконечное пополнение групп дает некоторую стабильность
взамен хаоса, свойственного дискретным группам. Мир симплициальных проконечных пространств, имея много свойств, схожих с дискретными, тем
не менее, крайне специфичен. Уже в размерности 1 (проконечные графы) стандартное утверждение о том, что в каждом связном графе существует
максимальное связное поддерево, уже не верно, см. \cite[Пример.~2.4]{Z}.

\section{Фильтрации.}

В этой и следующей главах работы мы сконцентрируем свое внимание на свободных симплициальных про-р-резольвентах про-р-групп.
Основным механизмом исследований станут $\Delta$-адическая фильтрация пополненного группового кольца и фильтрация Цассенхауза. Кратко
напомним ее основные конструкции, которые можно найти в \cite{L} и удачном обзоре, содержащемся в диссертации \cite{Gar}, а также введём
необходимые понятия, подчёркивающие ряд свойств пополненной групповой алгебры и её градуировки, см. также \cite{Qui5}.
\begin{zamech}
 Многие результаты этой и следующей глав имеют место и для фильтраций про-р-групп с помощью нижнего центрального ряда \cite{L}, однако, с
 нашей точки зрения, mod-p-вычисления в случае про-р-групп несколько удобней и информативней.
\end{zamech}

\begin{opr}
Пусть $G$ - про-р-группа, $\Delta^n(G)$ - n-ая степень аугментационного идеала её пополненной групповой алгебры $\mathbb{F}_p[[G]]$ и
$\nu$-фильтрация $\mathbb{F}_p[[G]]$, заданная убывающей цепочкой идеалов $(\Delta^n(G))_{n \in \mathbb{N}};$ т.е. $\nu(x) = \sup_{n: x
\in \Delta^n(G)} n$ для $x \in \mathbb{F}_p[[G]]$.
Фильтрацией Цассенхауса про-р-группы $G$ называется система замкнутых подгрупп в $G$, определенных по формуле $\mathcal{M}_n(G) = \{g \in
G\,|\, g-1 \in \Delta^n(G)\}$
\end{opr}

\begin{predl}\cite{KH}
Если $G$ - конечнопорожденная про-р-группа, то подгруппы $\mathcal{M}_n(G)$ открыты и определяют базис $1\in G$.
\end{predl}

Группы $\mathcal{M}_n(G)$ можно, следуя Дженнингсу, определить рекурсивно
\begin{predl}
Пусть $G$ - про-р-группа, тогда подгруппы $\mathcal{M}_n(G)$ могут быть описаны рекурсивно следующей формулой
$\mathcal{M}_n(G) = \overline{\mathcal{M}_{[n/p]}^p \prod_{i+j=n}[\mathcal{M}_i(G), \mathcal{M}_j(G)]}, \, n \in \mathbb{N},$
где $[n/p]$ - наименьшее целое число $\ge n/p.$
\end{predl}

Предложение 16 влечёт

(i) Если для каждого $i \ge 0$ обозначить $gr_i(G) =  \mathcal{M}_i(G)/\mathcal{M}_{i+1}(G)$, то $gr_i(G)$ - абелевы про-р-группы
экспоненты p>0, т.е. проконечные векторными $\mathbb{F}_p$-пространства.

(ii) Положим $\widehat{gr} G = \prod_{i \ge 0} gr_i G.$ Отображения $G \times G \to G,\, (x,y) \mapsto [x,y]$ и $G \to G,\, x \mapsto
x^p$
индуцируют корректно определенные отображения

$$[\cdot,\cdot]: gr_i G \times gr_j G \to gr_{i+j} G, \quad \cdot^{[p]}: gr_i G \to gr_{pi} G.$$
Распространим эти определения по линейности на $\widehat{gr} G:$
$$[\cdot,\cdot]: \widehat{gr} G \times \widehat{gr} G \to \widehat{gr} G, \quad
\cdot^{[p]}: \widehat{gr} G \to \widehat{gr} G.$$
Помимо фильтрации Цассенхауса нам потребуется $\Delta$-адическая фильтрация пополненного группового кольца про-р-группы и ассоциированный
с ней проконечный градуированный объект
$\widehat{gr} \mathbb{F}_p[[G]] = \prod_{n \ge 0} \frac{\Delta^n (G)}{\Delta^{n+1} (G)},$
который назовём проконечным ассоциированным групповым кольцом пополненного группового кольца про-р-группы $G$.
Будем обозначать через $\rho_n$ естественную сюръекцию $\rho_n: \mathcal{M}_n(G) \rightarrow \mathcal{M}_n(G)/ \mathcal{M}_{n+1}(G) =
gr_n G.$
Из определения фильтрации Цассенхауса группы $G$ следует существование отображения $\widehat{gr} G \to \widehat{gr} \mathbb{F}_p[[G]],$
заданного по формуле $\rho_n x \mapsto \rho_n(x-1).$

Введенные выше понятия можно несколько обобщить. Пусть $R$ - аугментированная проконечная алгебра над полем $\mathbb{F}_p$ с
аугментационным идеалом $\overline{R}$. Под фильтрацией алгебры $R$ будем подразумевать убывающую последовательность $R=F_0 R \supset F_1
R \supset...$ замкнутых подпространств таких, что $1\in R$ и $F_p R \cdot F_q R\subset F_{p+q} R.$ В этом случае $F_n R$ - двусторонний
идеал в $R$ и обозначим $\widehat{gr}R=\prod_{n=0}^{\infty} F_n/F_{n+1}$ - градуированный объект, который имеет естественную структуру
проконечной алгебры.

\begin{opr}
 Под полной аугментированной проконечной алгеброй мы будем подразумевать аугментированную проконечную алгебру $R$, наделенную фильтрацией
 $\{ F_n R \}$:

 (1) $F_1 R=\overline{R},$

 (2) $\widehat{gr}R$ порождается, как топологическая алгебра проконечным пространством $gr_1 R$,

 (3) $R= \varprojlim R/F_n R.$
\end{opr}

\underline{Класс полных аугментированных проконечных алгебр} образует категорию (PAA), где отображение $f:R\rightarrow R'$ - это
отображение аугментированных алгебр такое, что $f(F_n R)\subset F_n R'.$ Заметим, следуя Квиллену \cite[стр.267]{Qui5}, что в случае,
если $gr_1 R$ имеет конечную размерность, то $F_n R=\overline{R}^n.$ И как следствие, $R$ является фактором кольца формальных степенных
рядов $P=\mathbb{F}_p\langle \langle X_1,..,X_n\rangle \rangle$ от $n$ некоммутирующих переменных. Тогда любая PAA ввиду того, что она
представляется в виде проективного предела конечнопорождённых аугментированных алгебр тоже обладает свойством $F_n R=\overline{R}^n$ и
является фактором алгебры формальных степенных рядов над проконечным пространством $gr_1 R= \varprojlim X_\lambda,$ где $X_\lambda$ -
конечномерные. Можно считать, что по определению $\mathbb{F}_p\langle\langle X\rangle\rangle= \varprojlim \mathbb{F}_p\langle\langle
X_\lambda\rangle \rangle.$

Следуя Квиллену \cite[стр.270]{Qui5}, введём понятие проконечной алгебры Хопфа
\begin{opr}
Проконечной алгеброй Хопфа (PHA) называется PAA-алгебра $A$ с "диагональным" $\mathbb{F}_p$-модульным гомоморфизмом $\Delta :A\rightarrow
A \widehat{\otimes} A$ \cite[стр.269]{Qui5}. При этом структура биалгебры является коассоциативной и кокоммутативной. И наконец, имеется
$\mathbb{F}_p$-модульный "антипод" $s:A\rightarrow A,$ т.е. $\times (s \widehat{\otimes} 1)\Delta=\times (1 \widehat{\otimes}
s)\Delta=\varepsilon ,$ где $\varepsilon$ - аугментация.
\end{opr}
\begin{zamech}
Сохраним также термин проконечная алгебра Хопфа для градуированного объекта, связанного с фильтрацией $\bigoplus gr_n A$, тем самым
подчёркивая, что $gr_nA$ являются проконечными векторными пространствами.
\end{zamech}
\begin{prim}
Основной мотив введения понятия PHA - это пополненная $\mathbb{F}_p$-групповая алгебра про-р-группы $G$, т.к. диагональный гомоморфизм
$\Delta: G \rightarrow G \times G \quad (\Delta: g \mapsto g \times g)$ индуцирует непрерывный гомоморфизм пополненных групповых алгебр
$\Delta: \mathbb{F}_p[[G]] \rightarrow \mathbb{F}_p[[G \times G]]= \mathbb{F}_p[[G]] \widehat{\otimes} \mathbb{F}_p[[G]],$ антипод $s: g
\mapsto g^{-1}.$
\end{prim}
Если $A$ - PHA и $J$ - замкнутый Хопфовый идеал, т.е. $\Delta J \subset A \widehat{\otimes} J + J \widehat{\otimes} A,$ то $A/J$ является
PHA с отображением $\Delta,$ индуцированным из $A,$ т.к. $(R/J) \widehat{\otimes} (R/J)=R \widehat{\otimes} R/(R \widehat{\otimes} J + J
\widehat{\otimes} R).$

Если $A$ и $A'$ - PHA, то $A \widehat{\otimes} A'$ - PHA с проекциями $pr_1 :A \widehat{\otimes} A' \rightarrow A , \, pr_2 :A
\widehat{\otimes} A' \rightarrow A'$ и индуцированной аугментацией. Тогда $(A \widehat{\otimes} A',pr_1,pr_2)$ - прямое произведение $A$
и $A'$ в категории (PHA).

Введём проконечный аналог ограниченной алгебры Ли (обычное определение в \cite{DJ})
\begin{opr}
Проконечная алгебра Ли $\mathcal{L}$ (PLA) (т.е. проконечное $\mathbb{F}_p$-векторное пространство с непрерывной скобкой Ли)
характеристики $p>0$ называется ограниченной алгеброй Ли, если в $\mathcal{L}$ задано непрерывное отображение $a \mapsto a^{[p]}$,
удовлетворяющее следующим условиям:

(1) $[\alpha a]^{[p]}= \alpha^p a^{[p]},$

(2) $(a+b)^{[p]}=a^{[p]}+b^{[p]}+\sum\limits_{i=1}^{p-1} s_i(a,b)$, где $is_i(a,b)$ - коэффициенты при $\lambda^{i-1}$ в
$\textbf{ad}(\lambda a+b))^{p-1},$

(3) $[ab^{[p]}]=a(\textbf{ad}b)^p.$

\end{opr}

Любая PAA имеет структуру PLA, если в качестве скобки Ли рассмотреть $[x,y]=xy-yx$, а $a^{[p]}=a^p$, тогда условие (2) отражает
выполнение тождества на коэффициенты $(a+b)^{[p]}$.

Остановимся подробнее на фильтрации Цассенхауза свободных про-р-групп. Пусть задана свободная про-р-группа, порождённая проконечным
пространством $X$, тогда из аналогичного результата для конечнопорождённых про-р-групп \cite{L} переходом к проективному пределу
получаем

\begin{theorem}

Пусть $F(X)$ - свободная про-р-группа на проконечном пространстве $X$, тогда

(i) $gr {\mathbb{F}_p} [[F(X)]] \cong \textbf{T}(F/F^p[F,F]),$ где через $\textbf{T}(M)=\bigoplus_{n \geq 0}M^{\widehat{\otimes }n}$
обозначается тензорная алгебра проконечного $\mathbb{F}_p$-пространства $M$, и пополненный вариант  $( \widehat{gr} {\mathbb{F}_p}
[[F(X)]] \cong \widehat{\textbf{T}}(F/F^p[F,F]))$.

(ii) Имеет место изоморфизм проконечных ограниченных алгебр Ли $gr F(X) \xrightarrow{\sim} \mathcal{L}_{res}(X)$  $( \widehat{gr} F(X)
\xrightarrow{\sim} \widehat{\mathcal{L}}_{res}(X)),$
  $\widehat{\mathcal{L}}_{res}(X)$ - свободная проконечная ограниченная алгебра Ли над полем $\mathbb{F}_p,$ порождённая $X$, т.е.
  проективный предел свободных проконечных ограниченных алгебра Ли, порождённых конечными множествами; при этом, если $\xi_i =
  \rho_1(x_i),$ то $\xi_i \mapsto x_i.$
  И коммутативная диаграмма ограниченных алгебр Ли и алгебр Хопфа, где верхнее левое вложение задано на образующих правилом $x_i \mapsto
  1-x_i$,  $\widehat{gr}(\mathbb{F}_p[[F]])=\prod_{n \geq 0} gr_n \mathbb{F}_p[[F]],$ $\widehat{gr}(F)=\prod_{n \geq 0} gr_n(F),$
$$\xymatrix{ \widehat{gr}(F) \ar@{^(->}[r] \ar[dr]^{\sim} & \widehat{gr} \mathbb{F}_p[[F]] \ar[r]^{\sim} & \widehat{\textbf{T}}\langle
X\rangle\\
& \widehat{\mathcal{L}}_{res}(X) \ar@{^(->}[r] & \widehat{\mathcal{U}}_{L_{res}(X) \ar[u]^{\sim}}
}$$
\label{24}
\end{theorem}

Заметим, что нижнее вложение - это теорема Пуанкаре-Биркгофа-Витта, а крайний правый вертикальный изоморфизм - это изоморфизмы Квиллена,
собранные в проективный предел \cite{Qui4}.

Пусть задана свободная симплициальная про-р-группа $X$. Тогда гомоморфизмы граней и вырождения в $X$ индуцируют корректно-определенные
гомоморфизмы граней и вырождения в $\frac{X}{\mathcal{M}_n(X)}, \frac{\mathcal{M}_n(X)}{\mathcal{M}_{n+1}(X)},$ тем самым определяя
структуры симплициальных групп. Как известно \cite{KH}, гомоморфизмы групп индуцируют гомоморфизмы соответствующих пополненных групповых
колец. Тогда можно построить симплициальное групповое кольцо с операторами граней и вырождения, индуцированными соответствующими
гомоморфизмами исходной симплициальной группы. В нашем распоряжении оказываются проконечные симплициальные $\mathbb{F}_p$-алгебры
$\mathbb{F}_p[[X]]/\Delta^n(X)$ и $\frac{\Delta^n(X)}{\Delta^{n+1}(X)}.$
 Если $X$ - свободная симплициальная про-р-резольвента про-р-группы $G$, то возникают аугментированные симплициальные объекты
$\frac{X}{\mathcal{M}_n(X)} \twoheadrightarrow \frac{G}{\mathcal{M}_n(G)} \text{ и } \mathbb{F}_p[[X]]/\Delta^n(X) \twoheadrightarrow
\mathbb{F}_p[[G]]/\Delta^n(G),$
а также абелевы симплициальные про-р-группы экспоненты p>0:\
$\frac{\mathcal{M}_n(X)}{\mathcal{M}_{n+1}(X)} \twoheadrightarrow \frac{\mathcal{M}_n(G)}{\mathcal{M}_{n+1}(G)} \text{ и }
\frac{\Delta^n(X)}{\Delta^{n+1}(X)} \twoheadrightarrow \frac{\Delta^n(G)}{\Delta^{n+1}(G)} $ и их проконечные градуированные объекты.

В работе \cite{Mel1} был рассмотрен вопрос о существовании сюръективного гомоморфизма из про-р-группы с соотношением, являющимся
произведением коммутаторов и степеней образующих, в свободные про-р-группы конечного ранга.

Пусть $G$ - конечнопорожденная про-р-группа; внутренним про-р-рангом $Ir_p(G)$ группы $G$ назовем максимум рангов свободных про-р-групп,
являющихся эпиморфными образами $G$. Для элемента $r$ свободной про-р-группы $F$ конечного ранга положим $Ir_p(r,F) = Ir_p(F/R)$, где $R$
- порожденный $r$ нормальный делитель $F$. В дальнейшем будем писать $Ir_p(r)$ вместо $Ir_p(r,F)$.

\begin{theorem}\cite{Mel1}
Пусть $F$ - свободная про-р-группа ранга n=2k+m со свободной базой $\{x_1,...,x_k,y_1,...,y_k,z_1,...z_m\};$
$t = [x_1,y_1]...[x_k,y_k] z_1^{\al_1}...z_m^{\al_m},$
где $k \ge 0, m \ge 0, 0 \ne \al_i \in p\mathbb{Z}_p$ для $i=1,...,m.$ Тогда $Ir_p(t)=[n/2].$
\label{23}
\end{theorem}

Используя простейшие свойства фильтрации Цассенхауса симплициальной про-р-группы и элементарные методы линейной алгебры, покажем, следуя
идеям \cite{Sta}, как симплициальные методы могут быть применены  к исследованию задачи существования эпиморфизма конечнопорожденной
про-р-группы $G$ на свободную про-р-группу заданного ранга, в частности, воспроизведём ряд результатов \cite{Mel1}.

Пусть $G$ - про-р-группа с одним соотношением $G=\{x_1,...,x_n | r=1 \}$; будем считать, что $r \notin F^p[F,F]$, т.к. в противном
случае, согласно \cite{Se1}, его можно выбрать в качестве образующего $G$. Обозначим $F(m) = F(y_1,...y_m)$. Предположим, что существует
сюръективный гомоморфизм $\varphi: G \twoheadrightarrow F(m)$, тогда в соответствии с теоремой \ref{22} этот гомоморфизм индуцирует
некоторый гомоморфизм $\tilde{\varphi}$ симплициальной резольвенты $X \twoheadrightarrow G$ в симплициальную резольвенту  $C(F)$
свободной про-р-группы $F(m)$, в качестве которой можно выбрать постоянную симплициальную группу $F(m)$.

\begin{equation}
\xymatrix{
{...} \ar@<1ex>[r]^(0.25){d_2^2} \ar@<-1ex>[r]_(0.25){d_0^2} \ar[r] & F(n)*<\rho> \ar@<1ex>[r]^(0.6){d_1^1} \ar@<-1ex>[r]_(0.6){d_0^1}
\ar@{->>}[d]^{\varphi_1} & F(n) \ar[r] \ar@{->>}[d]^{\varphi_0} \ar[r] & G \ar[r] \ar@{->>}[d]^{\varphi} & 1\\
{...} \ar@<1ex>[r]^{id} \ar@<-1ex>[r]_{id} \ar[r]  & F(m) \ar@<1ex>[r]^{id} \ar@<-1ex>[r]_{id} & F(m) \ar[r] & F(m) \ar[r] & 1
}\label{15}
\end{equation}

Так как $F(m)$ - свободная про-р-группа, то можно построить ретракцию $\psi: F \to G$, т.е. такой гомоморфизм $\psi$ свободной
про-р-группы $F$ в $G$, что $\varphi \circ \psi = id_F$. Т.к. $X$ - асферичная симплициальная про-р-группа, то $\psi$ продолжается до
гомоморфизма симплициальных про-р-групп $\tilde{\psi}: C(F) \to X$. Откуда по теореме \ref{22} получаем некоторый гомоморфизм
$\tilde{\varphi} \circ \tilde{\psi}: C(F) \to C(F)$, который гомотопен тождественному гомоморфизму $C(F)$, т.е. $\tilde{\varphi} \circ
\tilde{\psi} = id_{C(F)}$.Но $C(F)$ - постоянная про-р-группа, поэтому композиция $\tilde{\varphi} \circ \tilde{\psi}$ просто совпадет с
тождественным гомоморфизмом, т.е. $\tilde{\varphi} \circ \tilde{\psi} \simeq id_{C(F)}$.
Т.к. $C(F)$ - постоянная про-р-группа, то $\tilde{\varphi}_1(\rho) = 1$ в \ref{15}.
Вложение свободной про-р-группы $F(n)$ в свое групповое кольцо $\mathbb{F}(n) \to \mathbb{F}_p[[F(n]]$ дает возможность определить
симплициальный гомоморфизм симплициальных групповых колец

\begin{equation}
\xymatrix{
{...} \ar@<1ex>[r]^(0.25){d_2^2} \ar@<-1ex>[r]_(0.25){d_0^2} \ar[r] & \mathbb{F}_p[[F(n)*\langle\rho\rangle]] \ar@<1ex>[r]^(0.6){d_1^1}
\ar@<-1ex>[r]_(0.6){d_0^1} \ar@{->>}[d]^{\tilde{\varphi}_1} & \mathbb{F}_p[[F(n)]] \ar[r] \ar@{->>}[d]^{\tilde{\varphi}_0} \ar[r] &
\mathbb{F}_p[[G]] \ar@{->>}[d]^{\varphi}\\
{...} \ar@<1ex>[r]^{id} \ar@<-1ex>[r]_{id} \ar[r]  & \mathbb{F}_p[[F(m)]] \ar@<1ex>[r]^{id} \ar@<-1ex>[r]_{id} & \mathbb{F}_p[[F(m)]]
\ar[r] & \mathbb{F}_p[[F]]
}\label{}
\end{equation}

Предположим, что $1-r \in \Delta^k / \Delta^{k+1},$ где $r = d_1^1(\rho),$ тогда
$1-r = \sum a_{i_1...i_k} (1-x_{i_1})...(1-x_{i_k})\mod \Delta^{k+1}.$
Рассмотрим симплициальные $\mathbb{F}_p$-векторные пространства $\frac{\Delta^n(X)}{\Delta^{n+1}(X)}$ и
$\frac{\Delta^n(C(F))}{\Delta^{n+1}(C(F))}$. Напомним \cite{KH}, что пополненное групповое кольцо $\mathbb{F}_p[[F]]$ свободной
про-р-группы может быть отождествлено с алгеброй формальных степенных рядов $\mathbb{F}_p(n)$ от $n$ независимых переменных (пополненная
тензорная алгебра $n$ - мерного $\mathbb{F}_p$-векторного пространства). Отождествление строится путем продолжения соответствия $x_i
\mapsto s_i - 1, i=1,...,n,$ где $s_i$ - образующие $F,$ до изоморфизма $\mathbb{F}_p(n) \simeq \mathbb{F}_p[[F]].$ Обозначим через $I$
аугментационный идеал $\mathbb{F}_p(n),$ который соответствует идеалу $\Delta$.

\begin{equation}
\xymatrix{
{...} \ar@<1ex>[r] \ar[r] \ar@<-1ex>[r] & \frac{I^k(\mathbb{F}_p(n+1))}{I^{k+1}(\mathbb{F}_p(n+1))} \ar@<1ex>[r]^(0.5){\bar{d}_1^1}
\ar@<-1ex>[r]_(0.5){\bar{d}_0^1} \ar@{->>}[d]^{\bar{\varphi}_1} & \frac{I^k(\mathbb{F}_p(n))}{I^{k+1}(\mathbb{F}_p(n))} \ar[r]
\ar@{->>}[d]^{\bar{\varphi}_0} \ar[r] & \frac{\mathfrak{g}^k}{\mathfrak{g}^{k+1}} \ar@{->>}[d]^{\varphi}\\
{...} \ar@<1ex>[r] \ar[r] \ar@<-1ex>[r]  & \frac{I^k(F(m))}{I^{k+1}(F(m))}  \ar@<1ex>[r]^{id} \ar@<-1ex>[r]_{id} &
\frac{I^k(F(m))}{I^{k+1}(F(m))} \ar[r] & \frac{f^k}{f^{k+1}}
}\label{}
\end{equation}

Будем отождествлять $I^k/I^{k+1}$ с векторным $\mathbb{F}_p$-пространством с базисом $\tilde{x}_{i_1}\otimes ...\otimes
\tilde{x}_{i_k}$.
Пусть $d_0^1(\rho) = 1$ и $d_1^1(\rho) = r,$ тогда в новых координатах $\bar{d}_1^1(1-\rho)=\sum a_{i_1...i_k} \tilde{x}_{i_1}\otimes ...
\otimes \tilde{x}_{i_k} = \eta(\tilde{x}_1,...,\tilde{x}_n).$
Т.к. $\tilde{\varphi}_1(1-\rho)=0$, то
\begin{equation}\eta (\tilde{\varphi}_1(\tilde{x}_1),...,\tilde{\varphi}_1(\tilde{x}_n))=0
\label{111}
\end{equation}

 Разложим $\tilde{\varphi}_1(\tilde{x}_i)$ в базисе $\tilde{y}_j \,$  $\tilde{\varphi}_1(\tilde{x}_i)=\sum_{j=1}^m b_{ij} \tilde{y}_j,\,
 i=1,...n$.
Существование ретракции $\psi$ позволяет утверждать, что $rank(B)=m,$ где $B=(b_{ij})$. Без ограничения общности можно считать, что $B$
имеет эшелонную форму. Этого можно добиться путем умножения $B$ слева на обратимую матрицу $C$ размера $m \times m$. Итак, необходимым
условием существования сюръективного гомоморфизма $G$ на свободную группу ранга $m$ является требование, что для формы $\eta$,
определенной в \ref{111}, существует матрица $B=(b_{ij})$ над $\mathbb{F}_p$ размера $n \times m$ и ранга $m$, что
$\eta(\sum_{j=1}^m b_{1j}\tilde{y}_j, ..., \sum_{j=1}^m b_{nj} \tilde{y}_j)=0.$
Т.е. $\sum_{i_1,...i_k}a_{i_1...i_k}b_{i_1j_1}...b_{i_kj_k} \tilde{y}_{j_1} \otimes ... \otimes \tilde{y}_{j_k}=0,$ т.е. для любого
набора $j_1,...j_k \in [1,m]$
\begin{equation}\sum_{i_1...i_k \in [1,n]} a_{i_1...i_k} b_{i_1j_1}...b_{i_kj_k} = 0\label{112}
\end{equation}
Пусть группа $G$ имеет представление
$P=\{x_1,...,x_n: x_1^e ... x_n^e = 1\},\,e=p^k.$

Учитывая, что характеристика $p > 0$, имеем разложение
$1-\tilde{r} = 1-(1-\tilde{x}_1)^e ...(1-\tilde{x}_n)^e = 1- (1-\tilde{x}_1^e)...(1-\tilde{x}_n^e)=(\tilde{x}_1^e+...+\tilde{x}_n^e) \mod
I^{e+1}.$
Таким образом, форма $\eta$ имеет все коэффициенты нулевые, кроме тех, в которых все индексы $i_1,...i_k$ одинаковые, в этом случае
коэффициенты равны единице.
Если $B=(b_{ij})$ - $n \times m$-матрица эшелонной формы, и мы возьмем набор из $k$ индексов
$j_1,...,j_k=s,s,s,...s,t$ в \ref{112}, то мы получим, что для всех $s,t \in [1,m]$
\begin{equation}\sum_{i=1}^n b_{is}^{e-1} b_{it} = 0.\label{113}
\end{equation}

Рассмотрим матрицу $C=(c_{si})$, в которой $c_{si} = b_{is}^{e-1},$ она получается из матрицы $B$ транспонированием, т.к. $B$ - эшелонная
матрица, то $C$ имеет тот же ранг $m$. Тогда равенство \ref{113} означает, что
\begin{equation}C \cdot B = 0.\label{114}
\end{equation}
Т.к. $dim_{F_p} \Ker C=n-m,$ а $dim_{\mathbb{F}_p}\Im B = m,$ то из \ref{114} следует, что $m \le n-m \Rightarrow m \le n/2.$

Можно отобразить нечетные образующие в самих себя, а четные образующие в обратные к предыдущим им нечетным (если последний окажется
нечетным, то он отобразится в единицу группы). Такое построение определяет ретракцию $G$ на свободную подгруппу ранга $(n-1)/2$ или $n/2$
в зависимости от того, какое из этих чисел целое. Т.е. $Ir_p(G) = [n/2]$.

Следуя \cite{Sta}, построим пример про-р-группы с одним соотношением, из которой не существует сюръективного гомоморфизма на свободную
про-р-группу ранга 2. Этот пример уже не покрывается теоремой \ref{23} и другими результатами \cite{Mel2}.

Для заданного $n$ пусть $\{a_{ij}\},\, 1 \le i < j \le n$ - это $\frac{n(n-1)}{2}$ попарно различных степеней двойки. Пусть $c$ - это
степень двойки, которая по крайней мере не меньше, чем каждая из $a_{ij}$. Определим $e_{ij}=c/a_{ij}.$ Рассмотрим соотношение
$r(x_1,...,x_n) = \prod_{i < j} [x_i^{a_{ij}}, x_j^{a_{ij}}]^{e_{ij}},$
мы предполагаем, что $x_1,..., x_n$ строго упорядочены. Пусть про-2-группа $G$ имеет представление
$G=\{x_1,...,x_n | r(x_1,...x_n) = 1\}.$
Как и в предыдущем примере рассмотрим свободную симплициальную резольвенту $G$ и соответствующий вариант с пополненными групповыми
кольцами. Заметим, что
$r(x_1,...x_n) = \sum_{i < j}(\tilde{x}_i^{a_{ij}} \tilde{x}_j^{a_{ij}}+\tilde{x}_j^{a_{ij}}\tilde{x}_i^{a_{ij}})^{e_{ij}} \mod
I^{2c+1}.$
Обозначим $\eta(\tilde{x}_1,...,\tilde{x}_n) = \sum_{i < j}(\tilde{x}_i^{\otimes a_{ij}} \otimes \tilde{x}_j^{\otimes
a_{ij}}+\tilde{x}_j^{\otimes a_{ij}} \otimes \tilde{x}_i^{\otimes a_{ij}})^{\otimes e_{ij}}.$

Если бы существовала ретракция $F(2) \to G \xrightarrow
{\varphi} F(2),$ мы бы получили  матрицу $B=(b_{ij})$ размера $n \times 2$, которая бы имела ранг 2, при этом
\begin{equation}
\eta(b_{11} \tilde{u} + b_{12} \tilde{v},...,b_{n1} \tilde{u} + b_{n2} \tilde{v}) = 0,
\label{115}
\end{equation}
где $u,v$ - образующие $F(2).$ Дальнейшие подробности и вычисления в \cite{Sta}.

\section{Спектральная последовательность.}

В своей работе \cite{Kan} Кан построил алгебраическую модель (симплициальные группы) для теории гомотопий, что позволило использовать
различные теоретико-групповые методы для дальнейших исследований. Кертис показал \cite{Cur71}, что фильтрация с помощью нижнего
центрального ряда приводит к спектральной последовательности, которая сходится к гомотопическим группам симплициальной группы. Спустя
короткое время Ректор \cite{Rec} рассмотрел фильтрацию Цассенхауза, которая привела к mod-p версии спектральной последовательности
Кертиса, сходящейся к факторам гомотопических групп симплициальной группы по модулю подгруппы элементов, порядки которых взаимнопросты с
p. Наконец, в работе \cite{Bou2} было доказано, что спектральная последовательность Кертиса-Ректора $\{E^i,X,d^ix \}$ в ряде случаев
совпадает при $i \ge 2 $ со спектральной последовательностью Адамса, см. также \cite{Sm}.

Построение спектральной последовательности производится следующим образом. Пусть $Y$ - односвязное симплициальное множество. Применяя
функтор Кана $Г$, получим свободную симплициальную группу $X=ГY.$ Фильтрации Цассенхауса и нижнего центрального ряда дают короткие точные
последовательности симплициальных групп
$$1 \to \mathcal{M}_{n+1}(X) \to \mathcal{M}_n(X) \to \mathcal{M}_n(X)/\mathcal{M}_{n+1}(X) \to 1, $$
$$1 \to \gamma_{n+1}(X) \to \gamma_n(X) \to \gamma_n(X)/\gamma_{n+1}(X) \to 1$$
Эти последовательности определяют длинные точные последовательности гомотопических групп, которые дают точные биградуированные пары. Как
мы уже отмечали, имеют место
\begin{theorem}~\cite{Cur71}
Пусть $K$ - связное, односвязное симплициальное множество, $G = ГK$ - его конструкция Кана, которая является по построению свободной
симплициальной группой. Тогда спектральная последовательность $E(G)$ сходится к $E^{\infty}(G) = \oplus_{r \ge 0} E_{r,q}^{\infty}$,
градуированной группе, ассоциированной с фильтрацией $\pi_q(ГK) = \pi_{q+1}(|K|).$ Группы $E^1(K)$ являются гомологическими инвариантами
$K$.
\end{theorem}
\begin{theorem}~\cite{Rec}
Пусть $K$ - связное и односвязное симплициальное множество, $G = ГK.$ Тогда $E^i(G,\mathbb{F}_p)$ слабо сходится к
$E^{\infty}(G,\mathbb{F}_p) = \oplus_{r \ge 0} E_{p,q}^{\infty}$, градуированной группе, ассоциированной с фильтрацией на $\pi_q(ГK) =
\pi_{q+1}(|K|)$  по модулю подгруппы элементов конечного порядка взаимнопростого с $p$.
\end{theorem}

Напомним стандартные построения \cite{Mak,MK}, связанные с точными парами.
\begin{opr}
Точная пара $\Upsilon$ - это пара $(D,E)$ модулей вместе с тремя морфизмами $i,j,k$ такими, что имеет место точный треугольник
$$\Upsilon: \xymatrix{D \ar[rr]^{i} & & D \ar[dl]^{j}\\
& E \ar[ul]^{k}}$$
\end{opr}

\begin{opr}(Производная пара).
Композиция $j \circ k: E \to E $ удовлетворяет тождеству $(j \circ k)(j \circ k) = j(k j) k = 0,$ т.е. мы можем образовать модуль $H(E) =
\frac{\Ker(j \circ k)}{\Im (jk)}.$

Построим треугольник
$$\Upsilon^1: \xymatrix{i(D) \ar[rr]^{i^1} & & i(D) \ar[dl]^{j^1}\\
& H(E) \ar[ul]^{k^1}},$$
\end{opr}
где $i^1$ - это ограничение $i$ на $i(D),$ а $j^1$ и $k^1$ заданы
$j^1(i(d)) = [j(d)], k^1([e]) = k(e).$
Заметим, что $j^1$ корректно определено, т.к. $i(d) = 0$ влечет, что для некоторого $e \in E$  существует $d = k(e)$ и $j(d) = j k(e)$ -
граница. Аналогично, $k^1$ корректно определено, т.к. если рассмотреть $e'' = e+jk(\tilde{e}),$ то $k(jk(\tilde{e})) = k j (k(\tilde{e}))
=0.$ Диаграмный поиск дает, что $\Upsilon^1$ тоже точная пара.

Итерация процесса получения точных пар приводит к $r$-ой производной паре $\Upsilon^r$ от $\Upsilon$
$$\Upsilon^r: \xymatrix{D^r \ar[rr]^{i^{(r)}} & & D^r \ar[dl]^{j^{(r)}}\\
& E^r \ar[ul]^{k^{(r)}}},$$
$D^r = i^r(D), E^r = H(E^{r-1}), \, i^{(r)}$ и $k^{(r)}$ индуцированы $i$ и $k$, так $j^{(r)}: [i^r(d)] \to [j(d)].$

Рассмотрим подробнее устройство дифференциала $d^1.$
По определению $d^1=j^1k^1=jk: E^1 \to E^1$ переводит $d$-гомологический класс $[y]$ в $d$-гомологический класс $[j(x)],$ где $x$ такой
элемент из $D$, что $i(x) = k(y)$ и $x$ существует, поскольку $y$ является циклом $j(k(y)) = 0.$ Дифференциал $d^2$ определен на
$E^2=H(E^1,d^1),$ причем $d^1[y] = 0,$ где $[y] \in E^1,$ тогда $d^2 = j^2 k^2$ переводит $d^1$-гомологический класс $[[y]] \in E^2$ в
$d^1$ - гомологический класс $[j^1(x)],$ где $x$ - такой элемент модуля $D^1$, что $i^1(x) = k^1[y].$
Таким образом, $d^2[[y]] = [j^1(i^1)^{-1} k^1[y]]=j(i)^{-1}k$  при условии, что первая часть формулы определена. Аналогично показывается,
что дифференциал $d^3$ можно рассматривать, как отображение, опеределенное там, где оба дифференциала $d^1$ и $d^2$ равны нулю, и
выражающееся формулой $d^3 = j^1(i^1)^{-1}(i^1)^{-1} k^1=j(i)^{-2}k.$

Дифференциал $d^n$ определен там, где дифференциалы $d,...d^{n-1}$ равны нулю, и выражается формулой

\begin{equation}
n \ge 1 \quad d^n = j(i)^{-(n-1)}k.
\label{10}
\end{equation}

\begin{predl}\cite{Mak}
$H(E) = \frac{k^{-1}(iD)}{j(\Ker(i))}$ и более общо $E^r = \frac{Z^r}{B^r},$ где $Z^r = k^{-1}(i^r D)$ и $B^r = j(\Ker (i^r)).$
\end{predl}

Пусть теперь $X$ - симплициальная про-р-группа. С этого места нас будет интересовать только mod-p вариант спектральной
последовательности. Рассмотрим точные последовательности симплициальных групп
$1 \to \mathcal{M}_n(X)/\mathcal{M}_{n+1}(X) \to X/\mathcal{M}_{n+1}(X) \to X/\mathcal{M}_{n}(X) \to 1,$
которые дают длинные точные последовательности гомотопических групп
$$...\to \pi_{i+1} (X/\mathcal{M}_n(X)) \to \pi_i(\mathcal{M}_n(X)/\mathcal{M}_{n+1}(X)) \to \pi_i(X/\mathcal{M}_{n+1}(X)) \to
\pi_i(X/\mathcal{M}_n(X)) \to ...$$

Такие длинные точные последовательности определяют, см. \cite{Hu,MK,Mak}, градуированные точные пары, которые несколько отличаются от
рассматриваемых в \cite{Rec,Cur71,Sm}
$$E_{n,m}^1 = \pi_m\hr{\frac{\mathcal{M}_n(X)}{\mathcal{M}_{n+1}(X)}} \quad D_{n,m}^1 = \pi_m\hr{\frac{X}{\mathcal{M}_n(X)}},$$
при этом
$$k_{n,m}^1: E_{n,m}^1 \to D_{n+1,m}^1 \quad \deg k^1 = (1,0),$$
$$i_{n+1,m}^1: D_{n+1,m}^1 \to D_{n,m}^1 \quad \deg i^1 = (-1,0),$$
$$j_{n,m}^1: D_{n,m}^1 \to E_{n,m-1}^1 \quad \deg j^1 = (0,-1).$$
Из \ref{10} ясно, что
$\deg d^r = (r,-1);  \, d^r: E_{n,m}^r(X) \to E_{n+r,m-1}^r(X), r \ge 1, $ но $\deg i^{(r)} = (-1,0),$ т.к. $i^{(r)}$ - это просто
ограничение $i$ на $i^{r-1}(D).$ Остается понять, какие степени $k^{(r)}$ и $j^{(r)}$, но $k^{(r)}$ просто совпадает с $k$, поэтому его
степень и будет степенью $k$, т.е. $\deg j^{(r)} = (1,0).$ Теперь мы можем вычислить $\deg j^{(r)}$, т.к. $\deg d^{(r)} = \deg
j^{(r)}+\deg k^{(r)} \Rightarrow (r,-1) = \deg j^{(r)}+(1,0) \Rightarrow \deg j^{(r)} = (r,-1)-(1,0)=(r-1,-1).$ В качестве примера
рассмотрим случай, когда $r = 2,$ тогда
$d_{n,m}^2: E_{n,m}^2 \to E_{n+2,m-1}^2;$
 $d_{n,m}^2$ по сути совпадает с $d^2 \approx j (i)^{-1} k,$ т.е.
$\pi_m \hr{\frac{\mathcal{M}_n(X)}{\mathcal{M}_{n+1}(X)}} \xrightarrow{k} \pi_m\hr{\frac{X}{\mathcal{M}_{n+1}(X)}} \xrightarrow{i^{-1}}
\pi_m\hr{\frac{X}{\mathcal{M}_{n+2}(X)}} \xrightarrow{j} \pi_{m-1}\hr{\frac{\mathcal{M}_{n+2}(X)}{\mathcal{M}_{n+3}(X)}}.$

Если рассмотреть фильтрацию свободной PHA $\mathbb{F}_p[[X]]$ степенями аугментационного идеала, то аналогичным образом возникают длинные
точные последовательности
$$...\to \pi_{i+1}(\mathbb{F}_p[[X]])/\Delta^n(X)) \to \pi_i(\Delta^n(X)/\Delta^{n+1}(X)) \to \pi_i(\mathbb{F}_p[[X]]/\Delta^{n+1}(X))
\to \pi_i(\mathbb{F}_p[[X]]/\Delta^n(X)) \to ...$$
Им соответствуют спектральные последовательности с дифференциалами бистепени $(r,-1)$
$$\bar{E}_{n,m}^1(X) = \pi_m(\Delta^n(X)/\Delta^{n+1}(X)), \quad \bar{d}^r: \bar{E}_{n,m}^r(X) \to \bar{E}_{n+r,m-1}^r(X) $$
Естественные гомоморфизмы $x\mapsto 1-x,\, x \in X_n, \, n \ge 0$ индуцируют морфизм спектральных последовательностей, см.
\cite[стр.166]{Gru}
$ \kappa_*: E(X) \to \bar{E}(X).$

Далее нас будет интересовать случай, когда $X$ - свободная симлициальная про-р-группа, например, свободная симплициальная резольвента
про-р-группы $G$.

\begin{theorem}
Пусть $X$ - свободная симплициальная про-р-группа. Обозначим $\overline{X}$ mod-p-абелианизацию $X$, тогда фильтрации Цассенхауза
свободной симплициальной про-p-группы и степенями аугментационного идеала её пополненного группового кольца приводят к спектральным
последовательностям проконечных векторных пространств:

(i)$E_{n,m}^1 = \pi_m(\mathcal{L}_n (\overline{X})), \bar{E}_{n,m}^1 = \pi_m(\textbf{T}_n (\overline{X})),$
в частности,
$E_{n,0}^1 = \mathcal{L}_n H_0(\overline{X}), \, \bar{E}_{n,0}^1=\textbf{T}_n H_0(\overline{X}); $ $\, E_{1,m} \cong \bar{E}_{1,m} \cong
H_m(\overline{X})$

(ii)индуцированный отображением $x \mapsto 1-x$ морфизм $\kappa_{n,0}^1: E_{n,0}^1(X) \to \bar{E}_{n,0}^1(X)\text{- инъективен.}$

(iii)$E_{n,0}^{\infty} = E_{n,0}^n = \mathcal{L}_n(\pi_0(X)),
\bar{E}_{n,0}^{\infty} = \bar{E}_{n,0}^n = gr_n \mathbb{F}_p[[\pi_0(X)]] $
\end{theorem}

\textbf{Доказательство.}

(i) описание $E_{n,m}^1,\bar{E}_{n,m}^1$ уже было получено в теореме  \ref{24} ввиду свободы $X$. Т.к.
$\frac{\Delta^n(X_m)}{\Delta^{n+1}(X_m)} \cong \prod^n \overline{X}_m,$ то проконечные варианты эквивалентностей Дольда-Кана (в том
смысле, что гомологии комплекса Мура симплициальной абелевой про-р-групп $\prod^n \overline{X}_m $ совпадают с гомологиями этой же
группы, рассматриваемой, как цепной комплекс с дифференциалами $d_m = \sum_{i=0}^m (-1)^i d_i^m $) и Эйленберга-Зильбера  \cite[$\S
29$]{May} (доказательства получаются разложением симплициальной про-р-группы в проективный предел конечных р-групп $\overline{X}=
\varprojlim \overline{X}_{\lambda} $, тогда $\prod^n \overline{X}_{\lambda} \simeq \overline{X}_{\lambda}^{\otimes n},$ откуда $\prod^n
\overline{X} \simeq \overline{X}^{\widehat{\otimes} n},$ вместе дают, что нормализованный комплекс абелевых про-р-групп
$...\xrightarrow{d_{m+1}} \frac{\Delta^n(X_m)}{\Delta^{n+1}(X_m)} \xrightarrow{d_m}  \frac{\Delta^n(X_{m-1})}{\Delta^{n+1}(X_{m-1})}
\xrightarrow{d_{m-1}} ...$
с дифференциалами $d_m = \sum_{i=0}^m (-1)^i d_i^m$ гомотопически эквивалентен $n$-кратному тензорному произведению комплексов
$\overline{X}$ $(\overline{X}^{\widehat{\otimes }n},\partial),$ где дифференциалы $\partial_m$ определяются индуктивно. Т.е., если
$\partial$ задан на $(n-1)$-тензорном произведении, то определим $\partial$ на $n$-тензорном произведении по формуле
$\partial (x \widehat{\otimes} y) = \partial(x) \widehat{\otimes }y + (-1)^{\deg x} x \widehat{\otimes} \partial(y),$
где $x \in \overline{X}, y \in  \overline{X}^{\widehat{\otimes} (n-1)}$. В частности, это позволяет использовать формулу Кюннета для
вычисления гомологий комплекса $\frac{\Delta^n(X)}{\Delta^{n+1}(X)}.$ Далее, $\pi_0\hr{\frac{\mathcal{M}_n(X)}{\mathcal{M}_{n+1}(X)}}
\cong \pi_0 (\mathcal{L}_n(\bar{X}))\cong \mathcal{L}_n(\pi_0(\bar{X})),$ т.к. функтор $\mathcal{L}$ сохраняет коуравнители. Аналогично,
т.к. функтор $\textbf{T}$ сохраняет коуравнители, то
$\pi_0\hr{\frac{\Delta^n(X)}{\Delta^{n+1}(X)}} \cong \pi_0(\textbf{T}_n(\bar{X})) \cong \textbf{T}_n(\pi_0(\bar{X})).$

(ii)
Представим $X\cong \varprojlim X_\lambda,$ где $X_\lambda$ - конечные симплициальные р-группы, тогда и $\mathbb{F}_p[X] \cong \varprojlim
\mathbb{F}_p[ X_\lambda].$ По теореме Квиллена \cite{Qui4} для каждой группы $X_\lambda$ имеет место изоморфизм $\mathcal{U}(grX_\lambda)
\cong gr\mathbb{F}_p[X_\lambda].$
 Ввиду теоремы Пуанкаре-Биркгофа-Витта имеем вложение $gr X_{\la} \hookrightarrow  \mathcal{U}(grX_\lambda)$. Т.к. функтор $\varprojlim$
 точен в категории проконечных пространств, то $gr X \hookrightarrow  gr\mathbb{F}_p[[X]],$ но функтор $\mathcal{U}$ сохраняет
 коуравнители, откуда и следует утверждение.

(iii)
$$\tiny{\xymatrix{ & &  {\pi_1\hr{\frac{X}{\mathcal{M}_{n-3}(X)}}} \ar@<-1ex>[d]_{i^{-1}} &
{\pi_1\hr{\frac{\mathcal{M}_{n-3}(X)}{\mathcal{M}_{n-2}(X)}}} \ar[dl]_{k_1} \ar[r]^{d^1} \ar@/_-5pc/[ddd]^{d^3} \ar[dr]^{d^2} &
{\pi_0\hr{\frac{\mathcal{M}_{n-2}(X)}{\mathcal{M}_{n-1}(X)}}}\\
& & {\pi_1\hr{\frac{X}{\mathcal{M}_{n-2}(X)}}}  \ar@<-1ex>[u]_i \ar[urr]^{j_1} \ar@<-1ex>[d]_{i^{-1}} &
{\pi_1\hr{\frac{\mathcal{M}_{n-2}(X)}{\mathcal{M}_{n-1}(X)}}}  \ar[r]^{d^1} \ar[dl]_{k_1} \ar@/_-3pc/[dd]^{d^2} &
{\pi_0\hr{\frac{\mathcal{M}_{n-1}(X)}{\mathcal{M}_{n}(X)}}}\\
& & {\pi_1\hr{\frac{X}{\mathcal{M}_{n-1}(X)}}} \ar[urr]_{j_1} \ar@<-1ex>[u]_i \ar@<-1ex>[d]_{i^{-1}} &
{\pi_1\hr{\frac{\mathcal{M}_{n-1}(X)}{\mathcal{M}_{n}(X)}}} \ar@{^(->}[dl]_{k_1} \ar[d]^{d^1}\\ {...} \ar[r] &
{\pi_1\hr{\frac{X}{\mathcal{M}_{n+1}(X)}}} \ar[r]^{i_1} & {\pi_1\hr{\frac{X}{\mathcal{M}_{n}(X)}}} \ar[r]^{j_1} \ar@<-1ex>[u]_i &
{\pi_0\hr{\frac{\mathcal{M}_{n}(X)}{\mathcal{M}_{n+1}(X)}}} \ar[r]^{k_0} &  {\pi_0\hr{\frac{X}{\mathcal{M}_{n+1}(X)}}} \ar[r]^{i_0} &
{\pi_0\hr{\frac{X}{\mathcal{M}_{n}(X)}}} \ar[r] & 0
}}$$

  Т.к. спектральная последовательность бирегулярна, то $E_{n,0}^{\infty} = E_{n,0}^n$, применяя предложение 17, имеем
$E_{n,0}^{\infty} = E_{n,0}^n = \frac{k^{-1}(\Im i^n)}{j\Ker i^n},$ $j(\Ker i^n) = j \hr{\pi_1\hr{\frac{X}{\mathcal{M}_n(X)}}}, \text{
т.к. } i^n: \pi_1\hr{\frac{X}{\mathcal{M}_n(X)}} \to \pi_1\hr{\frac{X}{X}} = 0,$
$k^{-1}(\Im i^n) = k^{-1}\hr{\pi_0\hr{\frac{X}{\mathcal{M}_{n+1}(X)}}} \cong \pi_0 \hr{\frac{\mathcal{M}_n(X)}{\mathcal{M}_{n+1}(X)}};$
присмотревшись к нижней строке диаграммы, получаем
$E_{n,0}^n = \frac{k^{-1}(\Im i^n)}{j (\Ker i^n)} \cong \frac{\pi_0\hr{\frac{\mathcal{M}_n(X)}{\mathcal{M}_{n+1}(X)}}}{j \hr{\pi_1 \hr{
\frac{X}{\mathcal{M}_n(X)}}}} \cong \Im k_0 \cong \Ker\hr{\pi_0\hr{\frac{X}{\mathcal{M}_{n+1}(X)}} \to
\pi_0\hr{\frac{X}{\mathcal{M}_n(X)}}}$
Функтор $Г_n: X \to \frac{X}{\mathcal{M}_n(X)}$ сохраняет коуравнители, следовательно
$\pi_0(Г_n(X) \cong Г_n(\pi_0(X))) \Rightarrow $
$\Ker\hr{\pi_0\hr{\frac{X}{\mathcal{M}_{n+1}(X)}} \to \pi_0\hr{\frac{X}{\mathcal{M}_n(X)}}} \cong \Ker(Г_{n+1}(\pi_0(X)) \to
Г_n(\pi_0(X))) \cong \frac{\mathcal{M}_n(\pi_0(X))}{\mathcal{M}_{n+1}(\pi_0(X))}.$
Аналогично, $\bar{E}_{n,0}^{\infty} = \Ker\hr{\pi_0\hr{\frac{\mathbb{F}_p[[X]]}{\Delta^{n+1}(X)}} \to
\pi_0\hr{\frac{\mathbb{F}_p[[X]]}{\Delta^n(X)}}} \cong \Ker \hr{\frac{\pi_0(\mathbb{F}_p[[X]])}{\Delta^{n+1} \pi_0(\mathbb{F}_p[[X]])}
\to \frac{\pi_0(\mathbb{F}_p[[X]])}{\Delta^n \pi_0(\mathbb{F}_p[[X]])}} \cong gr_n \pi_0(\mathbb{F}_p[[X]]),$ но функтор пополненной
групповой алгебры сохраняет коуравнители, откуда $gr_n \pi_0(\mathbb{F}_p[[X]]) = gr_n \mathbb{F}_p[[ \pi_0(X)]]. \boxtimes$

Для вычисления $\bar{E}_{n,m}^1 $ может быть применена формула Кюннета
$\overline{E}_{n,m}^1 \cong \sum_{i+j =m} \bar{E}_{n-1,j}^1 \widehat{\otimes} \bar{E}_{1,j}^1.$
Откуда, например, $\bar{E}_{2,1}^1 \cong \pi_1\hr{\frac{ \Delta^2(X)}{ \Delta^3(X)}} \cong \pi_1(\overline{X} \widehat{\otimes}
\overline{X}) \cong \sum_{i+j=1} \pi_i(\overline{X}) \widehat{\otimes} \pi_j(\overline{X}) \cong (\pi_0(\overline{X}) \widehat{\otimes}
\pi_1(\overline{X})) \oplus (\pi_1(\overline{X}) \widehat{\otimes } \pi_0(\overline{X}))$.
C учетом результатов главы о производных функторах
{\small\begin{align}
\bar{E}_{n,m}^1(G) = \bigoplus_{i+j=m}\bar{E}_{n-1,i}^1\widehat{\otimes}\bar{E}_{1,j}^1 =\bigoplus_{i_1+...+i_n = m}
\bar{E}_{1,i_1}^1\widehat{\otimes} ... \widehat{\otimes} \bar{E}_{1,i_n}^1 = \bigoplus_{i_1+...+i_n=m}H_{i_1+1}(G) \widehat{\otimes} ...
\widehat{\otimes} H_{i_n+1}(G)
\label{11}
\end{align}}

Пусть $X$ - свободная симплициальная про-р-резольвента про-р-группы $G$. Т.к. $\widehat{Г}^p \overline{W}G$ тоже свободная симплициальная
про-р-резольвента про-р-группы $G$, с учётом симплициальной гомотопической эквивалентности резольвент получаем изоморфизм соответствующих
спектральных последовательностей. Следовательно, для изучения свойств $\overline{E}(X)$ достаточно рассматривать спектральную
последовательность для $\overline{E}(\widehat{Г}^p \overline{W}G)$. Разложив $\overline{W}G= \varprojlim \overline{W}G_\lambda$ в
проективный предел конечных симплициальных пространств, получаем сходимость спектральной последовательности для соответствующих свободных
дискретных симплициальных групп $Г\overline{W}G_\lambda$ \cite{Bou1}, \cite[\S~9]{Cur71}, \cite[ch.13]{We}.  Переход к про-р-пополнению
ввиду конечнопорождённости не меняет спектральную последовательность. Аналогичным образом можно работать со спектральной
последовательностью $E(\widehat{Г}^p \overline{W}G)$ для любой симплициальной про-р-группы $G$. Переходя к проективному пределу по
$\lambda$, получаем следующее утверждение

\begin{predl}[Сходимость и самопожирание]
Пусть $X$ - свободная симплициальная про-р-резольвента про-р-группы $G$, тогда  $E_{*,m}^1(X) \Rightarrow \pi_m(\overline{G})(=0), \,
\overline{E}_{*,m}^1(X) \Rightarrow H_m(X)\,(=0), \, m > 0.$
\end{predl}

Пусть $G$ - симплициальная про-р-группа,  $ gr \mathbb{F}_p[[G_i]]$ является PHA для $i \geq 0 $, тогда $gr
\mathbb{F}_p[[G]]=\bigoplus_{i\geq 0} gr \mathbb{F}_p[[G_i]]$ c дифференциалом $d^1_n= \sum^n_{i\geq 0} (-1)^i d^n_i$ является
дифференциальной PHA  в смысле \cite[\S 3]{Bro}. Т.к. дифференциал уважает фильтрацию, получаем спектральную последовательность алгебр
Хопфа \cite[\S 3]{Bro}.
\begin{predl}[Структура дифференциальной алгебры Хопфа]
Пусть $X$- симплициальная про-р-группа, тогда $\bar{E}^r(X),$ $r \ge 1$ - проконечные алгебры Хопфа, структура алгебры Хопфа на
$\bar{E}^{r+1}(X)$ получается из $\bar{E}^r(X).$
\end{predl}

Пусть $G$ - симплициальная про-р-группа, тогда $G$ представима в виде проективного предела конечных р-групп $G \cong \varprojlim
G_\lambda ,$ откуда $\overline{W}G \cong \varprojlim \overline{W}G_\lambda .$ Как и ранее, $\widehat{Г}^p \overline{W}G \cong \varprojlim
\widehat{Г}^p \overline{W}G_\lambda. $ Для каждого $\lambda$ имеет место коммутативная диаграмма:

$\xymatrix{
& &  \pi_q(\overline{W}G_{\la}) \ar[r]^{h_q} \ar[d]^{\approx} & H_q(\overline{W}G_{\la}, \mathbb{F}_p) \ar[d]^{\approx} \\
... \ar[r] & \pi_{q-1}(\mathcal{M}_2 (\widehat{Г}^p \overline{W}G_{\la})) \ar[r] & \pi_{q-1}(\widehat{Г}^p \overline{W}G_{\la})
\ar[r]^{P} & \pi_{q-1}\hr{\frac{\widehat{Г}^p \overline{W}G_{\la}}{\mathcal{M}_2(\widehat{Г}^p \overline{W}G_{\la})}} \ar[r] & ...}
$

Переходя к проективному пределу по $\lambda$, получаем $H_q(\overline{W}G,\mathbb{F}_p)\cong \pi_{q-1} \hr{\frac{\widehat{Г}^p
\overline{W}G}{\mathcal{M}_2(\widehat{Г}^p \overline{W}G)}}.$
\begin{predl}[Естественность дифференциала ] \cite[Теорема 9.4]{Cur71}

Пусть $G$ - симплициальная про-р-группа, $\bar{E}_{1,*}^1(X)$ - спектральная последовательность её кофибрантной замены. Пусть для
некоторого $x \in \bar{E}_{1,*}^1(X)$ имеем $d^1(x)=\sum x_i \widehat{\otimes} x_j.$  Тогда дифференциал совпадает с коумножением в
гомологиях классифицирующего пространства $\overline{W}G$, т.е. $\Delta_{ \ast} ( \bar{x})=\sum \overline{x}_i \widehat{\otimes} (-1)^j
\overline{x}_j,$ где $ \bar{x}, \bar{x_i}, \bar{x_j}\in H_q(\overline{W}G,\mathbb{F}_p)$ соответствуют элементам $x,x_i, x_j \in
\pi_{q-1} \hr{\frac{\widehat{Г}^p \overline{W}G}{\mathcal{M}_2(\widehat{Г}^p \overline{W}G)}} $ с помощью изоморфизма предыдущей
диаграммы. Напомним, что коумножение в гомологиях $\Delta_{\ast}: H_*(\overline{W}G) \to H_*(\overline{W}G) \widehat{\otimes}
H_*(\overline{W}G)$ определяется диагональю $\Delta: \overline{W}G \to \overline{W}G \times \overline{W}G$ и может быть выражено явно с
помощью формулы Александера-Уитни \cite[гл.VIII]{Mak}.
Это коумножение  коассоциативно и кокоммутативно.
\end{predl}
Предложение 20 позволяет отождествить спектральную последовательность $\overline{E}$ c Cobar-спектральной последовательностью
\cite{Rec1}. В частности, $\overline{E}^2$ может быть представлен с помощью дифференциальной гомологической алгебры в виде
$\overline{E}^2 \approx Cotor^{H_{ \ast}(\overline{W}G)}(\mathbb{F}_p,\mathbb{F}_p)$ и $ \{ E^r \}$ имеет структуру спектральной
последовательности модулей над алгеброй Стинрода.

\begin{zamech}
Пусть $G$ - конечная р-группа, тогда из следствия 3 получаем, что у $G$ есть свободная конечнопорождённая про-р-резольвента конечного
типа, например свободная симплициальная про-р-группа $\widehat{Г}^p\overline{W}G$. $\overline{E}(\widehat{Г}^p\overline{W}G)$ будет
спектральной последовательностью конечномерных векторных пространств, причём $\overline{E}_{*,m}^1 $ сходится к $ 0$ при $m > 0$, и в
процессе своей "жизнедеятельности" оставляет информацию на горизонтальной прямой $m=0$ в виде размерностных факторов группового кольца
исходной р-группы $\overline{E}^{\infty}_{\ast,0}=gr\mathbb{F}_p[G]$. Тут проявляется  асимптотическая связь между группами гомологий,
дифференциалами спектральной последовательности и условием конечности $G$, т.е. обрывом цепочки $\Delta^n(G)/ \Delta^{n+1}(G).$.
Взглянем на результаты  Голода-Шафаревича \cite{GSh}, Винберга \cite{V}, Koxa \cite{KH} с наших позиций. Напомним, что в работе
\cite{GSh} было получено необходимое условие  конечности про-р-группы, заданной $d(G)$ образующими и $r(G)$ соотношениями $r(G) >
\frac1{4}d(G)^2$. Развивая неоднородные методы Винберга \cite{V}, Кох \cite{KH} получил неравенство
$$1 \le -dc_{n-1} + \sum_{\nu=0}^n c_{\nu} r_{n-\nu},$$ где $d$ - количество образующих группы; $r_n$ - количество соотношений, лежащих в
$F_m/F_{m-1},$ где $F_m$ - фильтрация Цассенхауза $F$ - свободной про-р-группы с $d$ образующими; $c_n = \sum_{\nu=0}^n b_{\nu},$ где
$b_n  = dim_{F_p} \frac{\Delta^n(G)}{\Delta^{n+1}(G)}, r_0=b_0=c_0=1.$
Из этого неравенства с помощью асимптотики рядов Пуанкаре удается получить замечательный результат \cite{KH}, связывающий границу $r(G)$
со "сложностью" соотношений в копредставлении конечной $p$-группы
\begin{equation}
r(G) > \frac{d(G)^m}{m^m}(m-1)^{m-1}.
\label{12}
\end{equation}
В монографии \cite{P1},посвященной проконечной алгебраической гомотопии, Тим Портер ставит "титаническую" задачу получения детальной
информации об относительном количестве числа образующих и соотношений в р-группах гомотопическими методами:

{\scriptsize
"Finally we would like to raise the tantalising problem of the link between
this cohomology theory and profinite presentations of (pro)finite groups. The
powerful results of Golod and Safarevic showed that the cohomology
theory of finite p-groups yields detailed information on the relative numbers
of generators and relations needed. Now one has the additional information, is
it possible to analyse this information yet further perhaps to give even better knowledge of the deficiency of a p-group?"}

 C учётом вышесказанного, задачу Портера можно "гомотопизировать", как вопрос об асимптотических свойствах Cobar-спектральной
 последовательности. Хотелось бы связать неравенство \ref{12} и скорость сходимости в духе теоремы Кертиса о связности \cite{Cur2,Qui4}.
\end{zamech}

\end{document}